\documentclass[a4paper,11pt]{article}
\usepackage[T1]{fontenc}
\usepackage{amssymb}
\usepackage{amsfonts, amsmath}
\usepackage{graphicx}

\newtheorem{Theo}{Theorem}

\newtheorem{Exmp}{Example}

\numberwithin{equation}{section}

\title{Isospectral deformations, the spectrum of Jacobi matrices,
infinite continued fraction and difference operators. Application
to dynamics on infinite dimensional systems}

\author{\textbf{A. Lesfari}
\\\emph{Department of Mathematics}
\\\emph{Faculty of Sciences}
\\\emph{University of Choua\"{i}b Doukkali}
\\\emph{B.P. 20, 24000 El Jadida, Morocco}.
\\\emph{E. mail : lesfariahmed@yahoo.fr.}}

\date{}
\begin{document}
\maketitle

\emph{Abstract}. This paper is devoted to the study of some
connections between coadjoint orbits in infinite dimensional Lie
algebras, isospectral deformations and linearization of dynamical
systems. We explain how results from deformation theory,
cohomology groups and algebraic geometry can be used to obtain
insight into the dynamics of integrable systems. Another part will
be dedicated to the study of infinite continued fraction,
orthogonal polynomials, the isospectral deformation of periodic
Jacobi matrices and general difference operators from an algebraic
geometrical point of view. Some connections with Cauchy-Stieltjes
transform of a suitable measure and Abelian integrals are given.
Finally the notion of algebraically completely integrable systems
is explained, techniques to solve such systems are presented and
some interesting cases appear as coverings of such dynamical
systems. These results are exemplified by several problems of
dynamical
systems of relevance in mathematical physics.\\
\emph{MSC(2010)}: 37K10, 32C35, 14B12, 14H40, 35S05, 11A55, 37K20.\\
\emph{Key words}: Integrable systems, Lax equations, cohomology
group, spectral curves, deformation, Jacobian varieties, periodic
Jacobi matrices, difference operators, continued fraction.

\section{Introduction}

The discovery towards the end of the 19th century by Poincar\'{e}
[51] that most nonlinear dynamical systems are not completely
integrable marked the end of a long and fruitful interaction
between Hamiltonian mechanics and algebraic geometry and the
interest in this subject decreased for more than half a century.
In fact many algebraic geometrical results such that elliptic and
hyperelliptic curves, Abelian integrals, Riemann surfaces, etc.,
have their origin in problems of mechanics. Fortunately the
discovery, 50 years ago, by Gardner, Greene, Kruskal and Miura
[15] that the Korteweg-de Vries (KdV) equation [30] :
$$
\frac{\partial u}{\partial t}-6u\frac{\partial u}{\partial
x}+\frac{\partial^{3}u}{\partial x^{3}}=0,\qquad u(x,0)=u(x), x\in
\mathbb{R}
$$
could be integrated by spectral methods has generated an enormous
number of new ideas in the area of Hamiltonian completely
integrable dynamical systems. The resolution of this problem has
led to unexpected connections between mechanics, spectral theory,
Lie groups, algebraic geometry and even differential geometry,
which have provided new insights into the old mechanical problems
of last centuries and many new ones as well. Lax [35] showed that
this equation is equivalent to the so-called Lax equation :
$$
\frac{dA}{dt}=[A,B]\equiv AB-BA,
$$
where $A$ (Sturm-Liouville) and $B$ are the differential operators
in $x$ :
$$A=-\frac{\partial ^{2}}{\partial x^{2}}+u,\qquad
B=4\frac{\partial ^{3}}{\partial x^{3}}-3\left(u\frac{\partial
}{\partial x}+\frac{\partial u}{\partial x}\right).$$ Lax equation
means that, under the time evolution of the system, the linear
operator $A(t)$ remains similar to $A(0)$. So the spectrum of $A$
is conserved, i.e. it undergoes an isospectral deformation. The
eigenvalues of $A$, viewed as functionals, represent the integrals
(constants of the motion) of the KdV equation. Thereafter, Mc
Kean-van Moerbeke [44], Dubrovin-Novikov [12] solved the periodic
problem for the KdV equation (for $x\in S^{1}$) in terms of a
linear motion on a real torus. This torus is the real part of the
Jacobi variety of a hyperelliptic curve with branch points defined
by the simple periodic and anti-periodic spectrum of $A$. Also the
motion is a straight line in the variables of the well known
Abel-Jacobi map. A parallel theory related to Jacobi matrices had
its origin in the periodic Toda problem [57] (discretized version
of the KdV equation). Krichever [33] generalized these ideas to
differential operators of any order, inspired by special examples
of Zaharov-Shabat [61], among which is the important
Kadomtsev-Petviashvili (KP) equation [24] :
$$\frac{3}{4}\frac{\partial ^{2}u}{\partial y^{2}}=
\frac{\partial}{\partial x}\left(\frac{\partial u}{\partial
t}-\frac{1}{4}\left( 6u\frac{\partial u}{\partial
x}+\frac{\partial^{3}u}{\partial x^{3}}\right)\right).
$$

Also this theory was generalized to difference operators of any
order by van Moerbeke and Mumford [60]. A one-to-one
correspondence was established between curves of a certain type
and classes of isospectral difference operators. They worked out a
systematic method which provides an algebraic map from the
invariants manifolds defined by the intersection of the constants
of the motion to the Jacobi variety of an algebraic curve
associated to Lax equation. Subsequently, these ideas led Mumford
[47] to show the absence of isospectral (here means that the
spectrum is given for all Floquet multipliers) deformations for
Laplace-like two-dimensional periodic difference operators by
relating the Picard variety with the class of such isospectral
operators and by showing that for a generic class of such
operators, the Picard variety is trivial. Therefore isospectral
flows appear only in the case of one-dimensional operators. I
shall not discuss here Mumford's result on the absence of
isospectral deformations for two-dimensional difference operators.
However given a dynamical Hamiltonian system, it remains often
hard to fit it into any of those general frameworks. But luckily,
most of the problems possess the much richer structure of the
so-called algebraic complete integrability (concept introduced et
systematized by Adler and van Moerbeke). A dynamical system is
algebraic completely integrable in the sense of Adler-van Moerbeke
[5, 7, 40] if it can be linearized on a complex algebraic torus
(Abelian variety).

Currently, the problem of finding and integrating nonlinear
dynamical Hamiltonian systems, has attracted a considerable amount
of attention in recent decades. Beside the fact that many such
systems have been on the subject of powerful and beautiful
theories of mathematics, another motivation for its study is : the
concepts of integrability which are applying to an increasing
number of physical systems, biological phenomena, population
dynamics, chemical rate equations, to mention only a few. However,
it seems still hopeless to describe or even to recognize with any
facility, those nonlinear systems which are integrable, though
they are exceptional. It is well known that the classical approach
to solving nonlinear integrable dynamical systems was dominated by
the question whether such systems can be solved by quadratures,
i.e., by a finite number of algebraic operations including the
inverting of functions. The method was based on solving the
Hamilton-Jacobi equation by separation of variables, after an
appropriate change of coordinates; for every problem finding this
transformation required a great deal of ingenuity. The solutions
of these problems can be expressed in terms of theta functions
related to Riemann surfaces. It must be emphasized that the
classical approach to proving that a system is integrable by
quadratures (in terms of hyperelliptic integrals) was something
very unsystematic and required a great deal of luck and ingenuity;
Jacobi himself was very much aware of this difficulty and in his
famous [22].

This paper is organized as follows : Section 1 is an introduction
to the subject. Section 2 concerns nonlinear integrable dynamical
systems which can be written as Lax equations with a spectral
parameter. Such equations have no a priori Hamiltonian content.
However, through the Adler-Kostant-Symes construction, we can
produce Hamiltonian dynamical systems on coadjoint orbits in the
dual space to a Lie algebra whose equations of motion take the Lax
form. We outline an algebraic-geometric interpretation of the
flows of these systems, which are shown to describe linear motion
on a complex torus. The relationship between spectral theory and
these systems is a fundamental aspect of the modern theory of
nonlinear integrable dynamical systems (see for example [7,39,
45]). We present a Lie algebra theoretical schema leading to
integrable systems, based on the Kostant-Kirillov coadjoint
action. Many problems on Kostant-Kirillov coadjoint orbits in
subalgebras of infinite dimensional Lie algebras (Kac-Moody Lie
algebras) yield large classes of extended Lax pairs. A general
statement leading to such situations is given by the
Adler-Kostant-Symes theorem [1, 31, 55] and the van
Moerbeke-Mumford linearization method provides an algebraic map
from the complex invariant manifolds of these systems to the
Jacobi variety (or some subabelian variety of it) of the spectral
curve. The complex flows generated by the constants of the motion
are straight line motions on these varieties. We present also the
Griffith's linearization method based on the observation that the
tangent space to any deformation lies in a suitable cohomology
group and that on algebraic curves, higher cohomology can always
be eliminated using duality theory. We explain how results from
deformation theory and algebraic geometry can be used to obtain
insight into the dynamics of integrable systems. These conditions
are cohomological and the Lax equations turn out to have a natural
cohomological interpretation. These results are exemplified by
several problems of dynamical integrable systems : Euler-Arnold
equations for the geodesic flow on the special orthogonal group
(the rotation group), Jacobi geodesic flow on the ellipsoid,
Neumann motion of a point on the sphere, Lagrange top, periodic
infinite band matrix, $n$-dimensional rigid body and Toda lattice.
The periodic Toda lattice consists of isospectral deformations of
periodic Jacobi matrices. The first flow describes a periodic
chain of particles interacting with an exponential potential. The
flow is conjugated to a motion of an auxiliary spectrum. Jacobi's
map transforms this motion into a linear flow on the Jacobi
variety of the hyperelliptic curve attached to the matrix. The
system is periodic or quasiperiodic. Section 3 is devoted to the
study of infinite continued fraction, orthogonal polynomials, the
isospectral deformation of periodic Jacobi matrices and general
difference operators from an algebraic geometrical point of view.
Complex Jacobi matrices play an important role in the study of
asymptotics and zero distribution of formal orthogonal
polynomials. Similarly, some connections with Cauchy-Stieltjes
transform of a suitable measure and Abelian integrals are given.
In Section 4 the notion of algebraically completely integrable
Hamiltonian systems is explained and techniques to solve such
systems are presented. Some important problems will be studied
namely the periodic $5$-particle Kac-van Moerbeke lattice and the
generalized periodic Toda systems. Also some interesting cases of
integrable systems for example the Ramani-Dorizzi-Grammaticos
(RDG) series of integrable potentials, a generalized
H\'{e}non-Heiles Hamiltonian, which will be studied in this
section, appear as coverings of algebraic completely integrable
systems. The manifolds invariant by the complex flows are
coverings of Abelian varieties and these systems are called
algebraic completely integrable in the generalized sense. The
later are completely integrable in the sense of Arnold-Liouville.
Also we shall see how some algebraic completely integrable systems
can be constructed from known algebraic completely integrable in
the generalized sense.

The concept of algebraic complete integrability is quite effective
in small dimensions and has the advantage to lead to global
results, unlike the existing criteria for real analytic
integrability, which, at this stage are perturbation results. In
fact, the overwhelming majority of dynamical systems, Hamiltonian
or not, are non-integrable and possess regimes of chaotic behavior
in phase space. I shall not discus here the weaker notion of
analytic integrability; the perturbation techniques developed in
that context are of a totally different nature. In recent years,
other important results have been obtained following studies on
the KP and KdV hierarchies. The use of tau functions related to
infinite dimensional Grassmannians, Fay identities, vertex
operators and the Hirota's bilinear formalism led to obtaining
important results concerning these algebras of infinite order
differential operators. In addition many problems related to
algebraic geometry, combinatorics, probabilities and quantum gauge
theory,..., have been solved explicitly by methods inspired by
techniques from the study of dynamical integrable systems. An
account of these results will appear elsewhere.

\section{Coadjoint orbits in Kac-Moody Lie algebras, isospectral
deformations and linearization}

A Lax equation (with parameter $h$) is given by a differential
equation of the form
\begin{equation}\label{eqn:euler}
\dot A(t,h)=[A(t,h),B(t,h)] \text{ or }[B(t,h),A(t,h)],\quad
\left(^.\equiv\frac{d}{dt}\right)
\end{equation}
where
$$A(t,h)=\sum_{k=l}^{m}A_{k}(t)h^{k},\qquad
B(t,h)=\sum_{k=l}^{m}B_{k}(t)h^{k},$$ are matrices or operators
and depend on a parameter $h$ (spectral parameter) whose
coefficients $A_{k}$ and $B_{k}$ are matrices in Lie algebras. The
pair $(A,B)$ is called Lax pair. The bracket $[,]$ is the usual
Lie bracket of matrices. The equation (2.1) establishes a link
between the Lie group theoretical and the algebraic geometric
approaches to complete integrability of dynamical systems. The
solution to (2.1) has the form $$A(t)=g(t)A(0)g(t)^{-1},$$ where
$g(t)$ is a matrix defined as $$\dot g(t)=-A(t)g(t).$$ We form the
polynomial $$P(z,h)=\det (A-zI),$$ where $z$ is another variable
and $I$ the $n\times n$ identity matrix. We define the curve
(spectral curve) $\mathcal{C},$ to be the normalization of the
complete algebraic curve whose affine equation is $P(z,h)=0$.

\begin{Theo}
A flow of the type (2.1) preserves the spectrum of $A$ for every
$h\in\mathbb{C}$ and its characteristic polynomial
$P(z,h)\equiv\det \left(A-zI\right)$. The latter defines an
algebraic curve
\begin{equation}\label{eqn:euler}
\mathcal{C}=\{(z,h) : P(z,h)=0\},
\end{equation}
for almost all $(z,h)\in\mathcal{C}$, which is time independent,
i.e., its coefficients $tr\left(A^{n}\right)$ are integrals of the
motion (equivalently, the matrix $A$ undergoes an isospectral
deformation).
\end{Theo}

The matrix $A-zI$, has a one-dimensional null-space, defining a
holomorphic line bundle on the curve $\mathcal{C}$. Whenever the
entries of the $A_k$ are moving in time, the curve $\mathcal{C}$
doest not move, inducing a motion on the set of line bundles. The
set of holomorphic line bundles on an algebraic curve forms a
group for the operation of tensoring $\otimes$ and the full set
with a given topological type is parametrized by the points of a
$g$-dimensional complex algebraic torus, where $g$ is the genus of
the curve. This torus that we note, $Jac(\mathcal{C})$, is the
Jacobian or Picard variety of the curve. When $\mathcal{C}$ is an
elliptic curve, $Jac(\mathcal{C})$ is isomorphic to $\mathcal{C}$.
Since the flow (2.1) induces deformations of line bundles, their
topological type remains unchanged and therefore it induces a
motion on the Jacobian variety; under some checkable condition on
$A$ and $B$, du to Griffiths [18] (see further for details). In
addition, some flows on Kostant-Kirillov coadjoint orbits in
subalgebras of infinite dimensional Lie algebras (Kac-Moody Lie
algebras) yield large classes of extended Lax pairs. A general
statement leading to such situations is given by the
Adler-Kostant-Symes theorem :

\begin{Theo}
Let $\mathcal{L}$ be a Lie algebra paired with itself via a
non-degenerate, ad-invariant bilinear form $\langle,\rangle$,
$\mathcal{L}$ having a vector space decomposition
$\mathcal{L}=\mathcal{K}+\mathcal{N}$ with $\mathcal{K}$ and
$\mathcal{N}$ Lie subalgebras. Then, with respect to
$\langle,\rangle$, we have the splitting
$\mathcal{L}=\mathcal{L}^{*}=\mathcal{K}^{\perp
}+\mathcal{N}^{\perp }$ and $\mathcal{K}^{\perp
}\approx\mathcal{N}^{*}$($\equiv$ the dual of $\mathcal{N}$)
paired with $\mathcal{N}$ via an induced form $\left\langle
\left\langle ,\right\rangle \right\rangle$ inherits the coadjoint
symplectic structure of Kostant and Kirillov; its Poisson bracket
between functions $H_{1}$ and $H_{2}$ on $\mathcal{N}^{*}$ reads
$$\left\{H_{1},H_{2}\right\}(a)=\left\langle \left\langle a,
\left[\nabla_{\mathcal{N}^{*}}H_{1},\nabla_{\mathcal{N}^{*}}H_{2}\right]
\right\rangle\right\rangle, \quad a\in \mathcal{N}^{*}.$$ Let
$V\subset \mathcal{N}^*$ be an invariant manifold under the above
coadjoint action of $\mathcal{N}$ on $\mathcal{N}^{*}$ and let
$\mathcal{A}(V)$ be the algebra of functions defined on a
neighborhood of $V$, invariant under the coadjoint action of
$\mathcal{L}$ (which is distinct from the $\mathcal{N-N}^*$
action). Then the functions $H$ in $\mathcal{A}(V)$ lead to
commuting vector fields of the Lax isospectral form, $$\dot
a=\left[a,pr_{\mathcal{K}}(\nabla H)\right],$$ $pr_{\mathcal{K}}$
projection onto $\mathcal{K}$.
\end{Theo}

The recent paper [8] gives the most general form of the
Adler-Kostant-Symes theorem. This theorem produces dynamical
Hamiltonian systems having many commuting integrals; some precise
results are known for interesting classes of orbits in both the
case of finite and infinite dimensional Lie algebras. The
finite-dimensional Lie algebras usually lead to noncompact
systems, and the infinite-dimensional ones to compact systems. Any
finite dimensional Lie algebra $\mathcal{L}$ with bracket $\left[
,\right]$ and Killing form $\left\langle ,\right\rangle$ leads to
an infinite dimensional formal Laurent series extension
$$\mathcal{L}=\sum_{-\infty }^{m}A_{i}h^{i}:A_{i}\in \mathcal{L},
\quad m\in \mathbb{Z} \mbox{ free},$$ with bracket $$\left[\sum
A_{i}h^{i},\sum
B_{j}h^{j}\right]=\sum_{i,j}\left[A_{i},B_{j}\right] h^{i+j},$$
and ad-invariant, symmetric forms $$\left\langle \sum
A_{i}h^{i},\sum B_{j}h^{j}\right\rangle_{k}=\sum_{i+j=-k}
\left\langle A_{i},B_{j}\right\rangle,$$ depending on $k\in
\mathbb{Z}$. The forms $\left\langle,\right\rangle_{k}$ are non
degenerate if $\left\langle,\right\rangle$ is so. Let
$\mathcal{L}_{p,q}$ $(p\leq q)$ be the vector subspace of
$\mathcal{L}$, corresponding to powers of $h$ between $p$ and $q$.
A first interesting class of problems is obtained by taking
$\mathcal{L}=\mathcal{G}l(n,\mathbb{R})$ and by putting the form
$\left\langle,\right\rangle_{1}$ on the Kac-Moody extension. Then
we have the decomposition into Lie subalgebras
$$\mathcal{L}=\mathcal{L}_{0,\infty }+\mathcal{L}_{-\infty
,-1}=\mathcal{K}+\mathcal{N},$$ with $\mathcal{K=K}^{\perp }$,
$\mathcal{N=N}^{\perp }$ and $\mathcal{K=N}^{*}$. Another class is
obtained by choosing any semi-simple Lie algebra $L$. Then the
Kac-Moody extension $\mathcal{L}$ equipped with the form
$\left\langle,\right\rangle=\left\langle,\right\rangle _{0}$ has
the natural level decomposition $$\mathcal{L}=\sum_{i\in
\mathbb{Z}}L_{i},\left[ L_{i,}L_{j}\right]\subset L_{i+j},
\quad\left[L_{0},L_{0}\right]=0, \quad L_{i}^{*}=L_{-i}.$$ Let
$B^{+}=\sum_{i\geq0}L_{i}$, $B^{-}=\sum_{i<0}L_{i}$. Then the
product Lie algebra $\mathcal{L\times L}$ has the following
bracket and pairing
\begin{eqnarray}
\left[\left(
l_{1},l_{2}\right),(l_{1}^{^{\prime}},l_{2}^{^{\prime}})\right]
&=&\left([l_{1},l_{1}^{^{\prime}}],-[l_{2},l_{2}^{^{\prime}}]\right),\nonumber\\
\left\langle\left(l_{1},l_{2}\right),(l_{1}^{^{\prime}},l_{2}^{^{\prime}})\right\rangle&=&\langle
l_{1},l_{1}^{^{\prime }}\rangle-\langle
l_{2},l_{2}^{^{\prime}}\rangle.\nonumber
\end{eqnarray}
It admits the decomposition into $\mathcal{K}+\mathcal{N}$ with
$$\mathcal{K}=\left\{ (l,-l):l\in \mathcal{L}\right\}, \quad
\mathcal{K}^{\perp }=\left\{ (l,l):l\in \mathcal{L}\right\},$$
$$\mathcal{N}=\left\{ (l_{-},l_{+}):l_{-}\in B^{-},l_{+}\in
B^{+},pr_{0}(l_{-})=pr_{0}(l_{+})\right\},$$$$ \mathcal{N}^{\perp
}=\left\{ (l_{-},l_{+}):l_{-}\in B^{-},l_{+}\in
B^{+},pr_{0}(l_{+}+l_{-})=0\right\},$$ where $pr_{0}$ denotes
projection onto $L_{0}$. Then from the last theorem , the orbits
in $\mathcal{N}^{*}\mathcal=K^{\perp }$ possesses a lot of
commuting Hamiltonian vector fields of Lax form. We state the
following theorem [2, 3, 60] :

\begin{Theo}
a) The invariant manifold $V_{m}$ $,$ $m\geq 1$ in
$\mathcal{K=N}^{*}$ , defined as
$$V_{m}=\left\{ A=\sum_{i=1}^{m-1}A_{i}h^{i}+\alpha h^{m}\text{ },
\text{ }\alpha=diag(\alpha_{1},\cdots ,\alpha _{n})\text{
fixed}\right\},$$ with $diag\left(A_{m-1}\right)=0$, has a natural
symplectic structure. The functions $$H=\left\langle
f(Ah^{-j}),h^{k}\right\rangle_{1},$$ on $V_{m}$ for good functions
$f$ lead to complete integrable commuting Hamiltonian systems of
the form
$$\dot{A}=\left[A,pr_{\mathcal{K}}(f^{\prime
}(Ah^{-j})h^{k-j})\right],$$ where
$$A=\sum_{i=0}^{m-1}A_{i}h^{i}+\alpha h^m,$$ and their trajectories
are straight line motions on the Jacobian of the curve
$\mathcal{C}$ of genus $\left( n-1\right)\left(nm-2\right)/2$
defined by (2.2). The coefficients of this polynomial provide the
orbit invariants of $V_{m}$ and an independent set of integrals of
the motion (of particular interest are the flows where $j=m,k=m+1$
which have the following form
\begin{equation}\label{eqn:euler}
\dot{A}=\left[A\text{ },\text{ }ad_{\beta \text{ }}ad_{\alpha
}^{-1}A_{m-1}+\beta h\right] ,\text{
}\beta_{i}=f^{\prime}\left(\alpha_{i}\right),
\end{equation}
the flow depends on $f$ through the relation
$\beta_{i}=f^{\prime}\left(\alpha_{i}\right)$ only).

b) (The van Moerbeke-Mumford linearization method) : The
N-invariant manifolds $$V_{-j,k}=\sum_{-j\leq i\leq
k}L_{i}\subseteq \mathcal{L\simeq K}^{\perp},$$ has a natural
symplectic structure and the functions $H(l_{1},l_{2})=f(l_{1})$
on $V_{-j,k}$ lead to commuting vector fields of the Lax form
$$\dot{l}=\left[ l,(pr^{+}-\frac{1}{2}pr_{0})\nabla H(l)\right],
\text{ }pr^{+}\mbox{ projection onto }B^{+},$$ their trajectories
are straight line motions on the Jacobian of a curve defined by
the characteristic polynomial of elements in $V_{-j,k}$ thought of
as functions of $h$, where $\nabla H(l)\in\mathcal{N}$ is the
gradient of $H$ thought of as a function on $\mathcal{L}$.
\end{Theo}

Using the van Moerbeke-Mumford linearization method [60], Adler
and van Moerbeke [3] showed that the linearized flow could be
realized on the Jacobian variety $Jac(\mathcal{C})$ (or some
subabelian variety of it) of the algebraic curve (spectral curve)
$\mathcal{C}$ associated to (2.1). We then construct an algebraic
map from the complex invariant manifolds of these dynamical
systems to the Jacobian variety $Jac(\mathcal{C})$ of the curve
$\mathcal{C}$. Therefore all the complex flows generated by the
constants of the motion are straight line motions on these
Jacobian varieties, i.e., the linearizing equations are given by
$$\int_{s_{1}(0)}^{s_{1}(t)}\omega _{k}+
\int_{s_{2}(0)}^{s_{2}(t)}\omega _{k}+\cdots +
\int_{s_{g}(0)}^{s_{g}(t)}\omega _{k}=c_{k}t\text{ },\text{ }0\leq
k\leq g,$$ where $\omega_{1},\ldots,\omega_{g}$ span the
$g$-dimensional space of holomorphic differentials on the curve
$\mathcal{C}$ of genus $g$.

\begin{Exmp} For $m=1$, i.e., for $V_1$, we choose
$$A=X+\alpha h, \quad X\in so(n).$$ In this case, the flow described by
equation (2.3) (where $\alpha_i$ and $\beta_i$ can be taken
arbitrarily) is reduced to the study of the Euler-Arnold equations
for the geodesic flow on $SO(n)$,
$$\dot{X}=[X,\lambda X],\quad (\lambda X)_{ij}=
\lambda_{ij}X_{ij},\quad
\lambda_{ij}=\frac{\beta_i-\beta_j}{\alpha_i-\alpha_j},$$ for a
left-invariant diagonal metric $\Sigma \lambda_{ij}X_{ij}$. The
natural phase space for this motion is an orbit defined in $SO(n)$
by $[n/2]$ orbit invariants. By Theorem 3, the problem is
completely integrable and the trajectories are straight lines on
$Jac(\mathcal{C})$ of dimension $(n-2)(n-1)/2$ and more
specifically, on the Prym variety
$Prym(\mathcal{C}/\mathcal{C}_0)\subset Jac(\mathcal{C})$ of
dimension $(n(n-1)/2-[w/2])/2$ induced by the natural involution
$\mathcal{C}\longrightarrow\mathcal{C}$, $(z, h)\longmapsto(-z,
-h)$, on $\mathcal{C}$ as a result of $X\in so(n)$;
$\mathcal{C}_0$ is the curve obtained by identifying $(z, h$) with
$(-z, -h)$.
\end{Exmp}

\begin{Exmp} For $m=2$, i.e., $V_2$, if one chooses
$$A=\alpha h^2-hx\wedge y-y\otimes y, \quad(x,y\in
\mathbb{R}^n),$$ which can also be considered as a rank 2
perturbation of the diagonal matrix $\alpha$ [45, 46, 2, 3], then
equation (2.3) reduces to
$$\dot{A}=[A,ad_\beta ad_\alpha^{-1}(y\wedge x)+\beta h],\quad\beta_i=f'(\alpha_i).$$
This equation can be reduced to the following nonlinear dynamical
system :
\begin{eqnarray}
\dot{x}&=&-(ad_\beta ad_\alpha^{-1}(y\wedge x)) x-\beta
y=-\frac{\partial H_\beta}{\partial
y},\nonumber\\
\dot{y}&=&-(ad_\beta ad_\alpha^{-1}(y\wedge x)) y=\frac{\partial
H_\beta}{\partial x},\nonumber
\end{eqnarray}
where
$$H_\beta=\frac{1}{2}\sum_i\beta_i\left(y_i^2+\sum_{j\neq
i}\frac{(x_iy_j-x_jy_i)^2}{\alpha_i-\alpha_j}\right),$$ which for
$f(z)=\ln z$, i.e., $\beta_i=\frac{1}{\alpha_i}$, we obtain the
problem of Jacobi geodesic flow on the ellipsoid :
$$\frac{x_1^2}{\alpha_1^2}+\cdots+\frac{x_n^2}{\alpha_n^2}=1,$$
expressing the motion of the tangent line ${x+sy :
s\in\mathbb{R}}$ to the ellipsoid in the direction $y$ of the
geodesic. For $f(z)=\frac{1}{2}z^2$, i.e., $\beta_i=\alpha_i$, we
get the Neumann motion [49] of a point on the sphere $S^{n-1}$,
$|x|=1$, under the influence of the force $-\alpha x$. From
theorem 3, both motions are straight lines on $Jac(\mathcal{C})$,
where $\mathcal{C}$ turns out to be hyperelliptic of genus $n-1$
(much lower than the generic one) ramified at the following $2n$
points : some point at $\infty$, the $n$ points $\alpha_i$ and
$n-1$ other points $\lambda_i$ of geometrical significance, based
on the observation that generically a line in $\mathbb{R}^n$
touches $n-1$ confocal quadrics. To be precise, the set of all
common tangent lines to $n-1$ confocal quadrics
$$Q_{\lambda_i}(x,x)+l=0, \quad i = 1, ..., n-1,$$ where
$Q_{u}(x,y)=\langle(u-\alpha)^{-1}x,y\rangle$, can be
parameterized by the quotient of the Jacobian of the hyperelliptic
curve $\mathcal{C}$ by an Abelian group $G$. The group is
generated by the discrete action obtained by flipping the signs of
$x_k$ and $y_k$ and some trivial one-dimensional action. Letting
$h\rightarrow0$ in the matrix $A$ and excising the largest
eigenvalue from this matrix leads to a new isospectral symmetric
matrix $$L=(I-P_y)(\alpha-x\otimes x)(I-P_y),$$ and a flow
$$\dot{L}=[ad_\beta ad_\alpha^{-1} x\wedge y,L],$$
where the spectrum of $L$ is given by the $n-1$ branch points
$\lambda_i$ above and zero. It follows that the tangent line
$\{x+sy : s\in\mathbb{R}\}$ to the ellipsoid remains tangent to
$n-2$ other confocal quadrics and the corresponding $n-1$
eigenfunctions of $L$ provide the orthogonal set of normals to the
$n-1$ quadrics at the points of tangency, hence recovering a
theorem of Chasles. The close relationship between Jacobi's and
Neumann's problems, which in fact live on the same orbits, was
implemented by Kn\"{o}rrer [29], who showed that the normal vector
to the ellipsoid moves according to the Neumann problem [49], when
the point moves according to the geodesic. These facts, as
investigated also by Kn\"{o}rrer [28] and others; the set of all
$n-1$ dimensional linear subspaces in the intersection of two
quadrics
\begin{eqnarray}
X_1^2+\cdots+X_n^2-Y_1^2-\cdots-Y_{n-1}^2&=&0,\nonumber\\
\alpha_1X_1^2+\cdots+\alpha_nX_n^2-\lambda_1Y_1^2-\cdots-\lambda_{n-1}Y_{n-1}^2&=&X_0^2,\nonumber
\end{eqnarray}
in $\mathbb{P}^{2n-1}$ is the Jacobian of the curve $\mathcal{C}$
defined above. This is done by observing that the set of linear
subspaces in the above quadrics is the same as the set of
$(n-2)$-dimensional linear subspaces tangent to $n-1$ quadrics
$$(\alpha_1-\lambda_j)X_1^2+\cdots+(\alpha_n-\lambda_j)X_n^2=X_0^2,\quad
j=1,2,...,n-1$$ which is dual to the set of tangents to the
confocal quadrics. The Neumann problem is also strikingly related
to the Korteweg-de Vries equation and various other nonlinear
dynamical systems (see [11]).
\end{Exmp}

\begin{Exmp} For another example of $V_2$ in theorem 3 , we consider the
Lagrange top (i.e., a symmetric top with a constant vertical
gravitational force acting on its center of mass and leaving the
base point of its body symmetry axis fixed) which evolves on an
orbit of type $V_2$; $n=3$, $$A=\Gamma+Mh+ch^2,$$ where $\Gamma\in
so(3)\simeq\mathbb{R}^3$ is the unit vector in the direction of
gravity and $M\in so(3)\simeq\mathbb{R}^3$ is the angular momentum
in body coordinates with regard to the fixed point; moreover
$c=(\lambda+\mu)\Upsilon$ where $\Upsilon\in
so(3)\simeq\mathbb{R}^3$ expresses the coordinates of the center
of mass and where $(\lambda+\mu, \lambda+\mu, 2\lambda)$ is the
inertia tensor in diagonalized form. The situation then leads to a
linear flow on an elliptic curve (see [54] and, for
higher-dimensional generalizations [53]).
\end{Exmp}

\begin{Exmp} As an example of $V_{ij}$ in theorem 3, b) (see [60, 2, 3]), we consider the
periodic infinite band matrix $M$ of period $n$ having $j+h+1$
diagonals; the spectrum of $M$ is defined by the points
$(z,h)\in\mathbb{C}^2$ such that $$Mv(h)=zv(h),\quad v(h)=(...,
h^{-1}v, v, hv,...), \quad v\in\mathbb{C}^n.$$ Let $M_h$ be the
square matrix obtained from $M$ and let $\mathcal{C}$ be the curve
defined by $\det(M_h-zI)=0$. Then the set of infinite band
matrices with $j+k+1$ diagonals, in higher dimensions many partial
results seem to lead to rigidity. In fact, it was shown that a
discrete 2-dimensional Laplacian cannot be deformed, given its
periodic spectrum; the proof can be summarized by the observation
that the Picard variety of most algebraic surfaces is trivial; the
proof that the specific spectral surface defined by the
$2$-dimensional Laplacian has trivial Picard variety is based on
the technique of toroidal embedding, which reduces cohomological
computations to combinatorial questions. Finally, inspired by the
dynamical systems, Mumford [48] has given a beautiful description
of hyperelliptic Jacobians of dimension $g$. Let $y^2=R(z)$ be the
monic polynomial of degree $2g+l$ defining the curve $\mathcal{C}$
and let $\theta$ be the  th\^{e}ta divisor. Then
$\mbox{Jac}(\mathcal{C})\backslash\theta$ is a variety of
polynomials $U$, $V$ with $\mbox{deg }U=g$, $\mbox{deg }V\leq g-1$
and $U$ monic such that $U|B-V^2$.
\end{Exmp}

As mentioned before, in an unifying approach Griffiths [18] has
found necessary and sufficient conditions on $B$ for the Lax flow
(2.1) to be linearizable on the Jacobi variety of its spectral
curve, without reference to Kac-Moody Lie algebras, so that the
flow of the Lax form (2.1) can be linearized on the Jacobian
variety $\mbox{Jac}(\mathcal{C})$ for $\mathcal{C}$ defined by
(2.2). Suppose that for every $p(z,h)$ belonging to the curve
$\mathcal{C}$, with $\dim\ker(A-zI)=1$, (i.e., the corresponding
eigenspace of $A$ is one-dimensional) and generated by a vector
$v(t,p)\in V$ where $V\simeq\mathbb{C}^n$ is an $n$-dimensional
vector space. There is then a family of holomorphic mappings which
send $(z,h)\in \mathcal{C}$ to $\ker (A-zI)$: $$f_t :
\mathcal{C}\longrightarrow \mathbb{P}V, \quad p\longmapsto
\mathbb{C}v(t,p),$$ called the eigenvector map associated to the
Lax equation. We set $$L_t
=f^*_t\left(\mathcal{O}_{\mathbb{P}V}(1)\right)\in
\mbox{Pic}^d(\mathcal{C})\cong\mbox{Jac}(\mathcal{C}), \quad
L=L_0,$$ where $d=\mbox{deg }f_t(\mathcal{C})$;
$\mathcal{O}_{\mathbb{P}V}(1)$ is the hyperplane line bundle on
$\mathbb{P}V$ and $\mbox{Pic}^d(\mathcal{C})$ the Picard variety
of $\mathcal{C}$, i.e., let us recall that it is the set of
straight bundles of degree $d$ on $\mathcal{C}$. By continuity,
the degree of $L_t$ does not vary with time $t$. Let $\textbf{H}$
be the hyperplane class of $\mathbb{P}V$. We have $$\mbox{deg
}L_t=\int_\mathcal{C}f_t^*\textbf{H}=\int_{f_t(\mathcal{C})}\textbf{H}.$$
This expression is the Poincar\'{e} dual of the class
$[\mathcal{C}]$ of $\mathcal{C}$ and coincides with the degree of
$\mathcal{C}$. Hence $\mbox{deg }L_t=\mbox{deg}(\mathcal{C})$.
While $t$ varies, $L_t$ moves in $\mbox{Pic}^d(\mathcal{C})$.
Therefore, if we fix a line bundle $L_0\in
\mbox{Pic}^d(\mathcal{C})$, the line bundle $L_0^{-1}\otimes L_t$
moves in the Jacobian variety
$$\mbox{Jac}(\mathcal{C})=H^1(\mathcal{C},\mathcal{O}_{\mathcal{C}})/H^1(\mathcal{C},\mathbb{Z})\simeq
H^0(\mathcal{C},\Omega_{\mathcal{C}})^*/H_1(\mathcal{C},\mathbb{Z}),$$
i.e., the mapping $L\longmapsto L_0^{-1}\otimes L$ induces a
morphism $\mbox{Pic}^d(\mathcal{C})\simeq
\mbox{Jac}(\mathcal{C})$. The motion of the line bundle
$L_0^{-1}\otimes L_t$ depends on the choice of the matrix $B$ and
a question arises : determine necessary and sufficient conditions
on the matrix $B$ so that the flow
\begin{equation}\label{eqn:euler}
t\longmapsto L_t\in \mbox{Jac}(\mathcal{C}),
\end{equation}
can be linearized on the Jacobian variety
$\mbox{Jac}(\mathcal{C})$. As we have pointed out, Griffiths has
found necessary and sufficient conditions of a cohomological
nature on $B$ that the flow $t\longmapsto L_t\in
\mbox{Jac}(\mathcal{C})$, be linear. His method is based on the
observation that the tangent space to any deformation lies in a
suitable cohomology group and that on algebraic curves, higher
cohomology can always be eliminated using duality theory. In fact,
applying more or less standard cohomological techniques from
deformation theory [9], we may give necessary and sufficient
conditions that the map $t\longmapsto L_t$ be linear.

Let
\begin{equation}\label{eqn:euler}
f:\mathcal{C}\longrightarrow X,
\end{equation}
be a non-constant holomorphic map where $\mathcal{C}$ is a given
smooth algebraic curve and $X$ is a complex manifold. We define
the normal sheaf of $\mathcal{C}$ in $X$ by the exact sequence
\begin{equation}\label{eqn:euler}
0\longrightarrow
\Theta_\mathcal{C}\overset{f_{\ast}}{\longrightarrow}f^*\Theta_X\longrightarrow
N_f\longrightarrow0
\end{equation}
with $\Theta_\mathcal{C}$, $\Theta_X$ are the respective tangent
sheaves and $f_{\ast}$ is the differential of $f$. Then the
Kodaira-Spencer tangent space [9] to the moduli space of the map
(2.5) is given by $H^0(\mathcal{C},N_f)$. If $$f_t :
\mathcal{C}\longrightarrow X, \quad f_0=f,$$ is a deformation of
(2.5), then $\dot{f}\in H^0(\mathcal{C},N_f)$ the corresponding
infinitesimal deformation at $t=0$, i.e., in local product
coordinates $(z,t)$ on $\cup_t\mathcal{C}_t$ and $w=(w^1,
w^2,...,w^n)$ of $X$, $f_t$ is given by $(t,\xi)\longmapsto
w(t,\xi)$, then the section $\dot{f}\in H^0(\mathcal{C},N_f)$ is
locally given by $\left.\frac{\partial w(t,\xi)}{\partial
t}\right|_{t=0} \mbox {modulo } \frac{\partial w(0,\xi)}{\partial
z}$. The corresponding cohomological sequence of (1.6) is
$$H^0(\Theta_\mathcal{C})\longrightarrow
H^0(f^*\Theta_X)\longrightarrow
H^0(N_f)\overset{\overline{\partial}}{\longrightarrow}H^1(\Theta_\mathcal{C}).$$
Here $H^1(\Theta_\mathcal{C})$ is the tangent space to the moduli
space of $\mathcal{C}$ as an abstract curve and
$$\overline{\partial}(\dot{f})\equiv\dot{\mathcal{C}}\in
H^1(\Theta_\mathcal{C}),$$ is the tangent to the family of curves
$\{\mathcal{C}_t\}$. Thus the tangent space to deformations of
(2.5) where the curve $\mathcal{C}$ remains fixed, is given by
$$H^0(f^*\Theta_X)/H^0(\Theta_\mathcal{C})\subset H^0(N_f).$$ Since
the isospectral curve $\mathcal{C}$ is independent of $t$, this is
the situation that we are interested in.

In the following, we take again the vector space $V$ of dimension
$n$ and assume that $X=\mathbb{P}V$ (projective space) and
consider the Euler sequence
$$
0\longrightarrow
\mathcal{O}_{\mathbb{P}V}\overset{i}{\longrightarrow}
V\otimes\mathcal{O}_{\mathbb{P}V}(1)\overset{p}{\longrightarrow}
\mathcal{O}_{\mathbb{P}V}\longrightarrow0
$$
This is an exact sequence of vector bundles, so that it remains
exact after pulling back to $\mathcal{C}$ via $f^*$ (combining
this with (2.6)). We have then a diagram of exact sequences
($L=f^*\mathcal{O}_{\mathbb{P}V}(1)$) :
$$\begin{array}{ccccccccc}
&&&0&&&&&\\
&&&\downarrow&&&&&\\
&&&\mathcal{O}_\mathcal{C}&&&&&\\
&&&\quad\left\downarrow v\right.&&&&&\\
&&&V\otimes L&&&&&\\
&&&\downarrow&&&&&\\
0&\longrightarrow&
\Theta_\mathcal{C}&\overset{f_{\ast}}{\longrightarrow}f^*\Theta_{\mathbb{P}V}&\longrightarrow&
N_f&\longrightarrow&0\\
&&&\downarrow&&&&&\\
&&&0&&&&
\end{array}
$$
The associated cohomology diagram contains the following piece :
$$\begin{array}{ccccccccc}
&&&H^0(\mathcal{C},V\otimes L)&&&&&\\
&&&\left\downarrow \tau\right.&&&&&\\
&&H^0(\mathcal{C},\Theta_\mathcal{C})&\longrightarrow
H^0(\mathcal{C},
f^*\Theta_{\mathbb{P}V})&\overset{j}{\longrightarrow}&H^0(\mathcal{C},N_f)&
\overset{\overline{\delta}}{\longrightarrow}&H^1(\mathcal{C},\Theta_\mathcal{C})\\
&&&\left\downarrow \delta\right.&&&&&\\
&&&H^1(\mathcal{C},\mathcal{O}_\mathcal{C})&&&&
\end{array}
$$
Consider the family of holomorphic maps $f_t :
\mathcal{C}\longrightarrow \mathbb{P}V$. Locally choose a
coordinate $\xi$ on $\mathcal{C}$ and a position vector mapping
$(t,\xi)\longmapsto v(t,\xi)\in V\backslash\{0\}$, i.e., a local
lift $v_t$ of $f_t$ to $V\backslash\{0\}$, such that
$$f_t(\xi)=\mathbb{C}.v(t,\xi)\subset V.$$ Notice that $v_t$ is a
time-dependent map $\mathcal{C}\longrightarrow V\backslash\{0\}$.
This lift is not canonical and exists only locally, but we are
going to use it to define an object denoted $\dot{v}$ which will
be independent of the lift and therefore will be globally well
defined. Since $\mathcal{O}_{\mathbb{P}V}$ is the tautological
bundle of $\mathbb{P}V$, the fibre of
$f^*\mathcal{O}_{\mathbb{P}V}(-1)$ at a point $p\in\mathcal{C}$
may be identified with the space $\mathbb{C}v_t(p)$, which defines
the maps
$f^*\mathcal{O}_{\mathbb{P}V}(-1)V\otimes\mathcal{O}_\mathcal{C}$
and $$v_t:\mathcal{O}_\mathcal{C}\longrightarrow V\otimes L_t,
\quad \phi\longmapsto\phi v_t,$$ where $v_0$ coincides with the
application $v$ mentioned in the previous diagram. If
$\widetilde{v}$ is another lift given by
$$\widetilde{v}(t,\xi)=\kappa(t,\xi)v(t,\xi), \quad \kappa\neq0,$$ then
we have $$\dot{\widetilde{v}}=\kappa\dot{v}+\dot{\kappa}v.$$ Set
$$\dot{v}(\xi)=\left.\frac{\partial v(t,\xi)}{\partial
t}\right|_{t=0} \mbox {modulo }v(t,\xi).$$ The latter quantity is
well defined of the representative position mapping of $v$, i.e.,
since the inclusion
$$\mathcal{O}_\mathcal{C}\overset{v}{\hookrightarrow}V\otimes L,
\quad L=f^*\mathcal{O}_{\mathbb{P}V}(1),$$ is locally given by
$\mathcal{O}_\mathcal{C}\ni\phi\longmapsto\phi.v$, it follows that
$$\dot{v}\in H^0(\mathcal{C},V\otimes
L/\mathcal{O}_\mathcal{C})=H^0(\mathcal{C},
f^*\Theta_{\mathbb{P}V}),$$ is a well-defined and independent of
the choice of the lift. Then we have $ j(\dot{v})=f$. We are
interested in the tangent vector
$$\dot{L}\equiv\left.\frac{d L_t}{d
t}\right|_{t=0}\in H^1(\mathcal{C},\mathcal{O}_\mathcal{C}).$$

\begin{Theo}
We have $$\dot{L}=\delta(\dot{v}),$$ where $\dot{v}$ is the
infinitesimal variation of
$f_t:\mathcal{C}\longrightarrow\mathbb{P}V$ and in particular,
$\dot{L}=0$ if and only if $\dot{v}=\tau(w)$ for some $w\in
H^0(V\otimes L)$ where $\tau$ is the map in diagram above.
\end{Theo}

We write
$$B(t,h)=\sum_{k=0}^{N}B_{k}(t)h^{k}=\sum_{k=0}^{N}B_{k}(t)h^{N-k}_0h^k_1,$$
where we have regarded $h$ as an affine coordinate on the
projective line $\mathbb{P}^1$ which is the base of the covering
$\pi : \mathcal{C}\longrightarrow\mathbb{P}^1$, while $h_0$, $h_1$
are homogeneous coordinates. Recall that $B(t,h)\in
H^0(\mathcal{C},\mbox{Hom}(V,V(N)))$ where $V$ is the sheaf of
sections of the trivial bundle $\mathcal{C}\times V$,
$V(D)=V\otimes\mathcal{O}_\mathcal{C}(D)$ (Here $B(t,h)$ is a
holomorphic section of the bundle
$\mbox{Hom}(V,V)\otimes\mathcal{O}_\mathcal{C}(N)$,
$\mathcal{O}_\mathcal{C}(N)=\pi^*\mathcal{O}_{\mathbb{P}^1}(N)$,
i.e., we are viewing $h=[h_0, h_1]$ as a homogeneous coordinate on
$\mathbb{P}^1$ pulled up to $\mathcal{C}$). Let $D=(h_0^N)$, be
the divisor $N.\pi^{-1}(\infty)$ on the curve $\mathcal{C}$.
Therefore $B/h_0^N\in H^0(\mathcal{C},\mbox{Hom}(V,V(D)))$, $v\in
H^0(V\otimes L)$ where $V(D)\cong V(N)$ are the sections of
$V\otimes\mathcal{O}_\mathcal{C}(D)$ (Here $B/h_0^N$ is a matrix
in $\mbox{Hom}(V, V)$ with meromorphic functions in
$H^0(\mathcal{C},\mathcal{O}_\mathcal{C}(D))$ as entries, i.e., we
are viewing $h_1/h_0$ as a function in
$H^O(\mathcal{C},\mathcal{O}_\mathcal{C}(D))$). Hence
$\left(\frac{B}{h_0^N}\right).v\in H^0(\mathcal{C},V\otimes L(D))$
and the cohomological interpretation of Lax equation is given by

\begin{Theo}
We have
$$\dot{v}=\tau\left(\displaystyle{\frac{B}{h_0^N}.v}\right),$$ and
$\dot{L}=0$, if and only if there is a meromorphic function
$\varphi\in H^O(\mathcal{C},\mathcal{O}_\mathcal{C}(D))$ such that
$\displaystyle{\frac{B}{h_0^N}.v+\varphi v}\in
H^0(\mathcal{C},V\otimes L(D))$ is holomorphic.
\end{Theo}

Near the point $p=(h,z)\in\mathcal{C}$, differentiating with
regard to $t$ the eigenvalue problem $Av(t,p)=zv(t,p)$, leads to
$\dot{A}v+A\dot{v}=z\dot{v}$. Using the Lax equation :
$\dot{A}=[B,A]$, we obtain $A(\dot{v}-Bv)=z(\dot{v}-Bv)$. Since
generically the eigenvalues have multiplicity $1$, we have
\begin{equation}\label{eqn:euler}
Bv=\dot{v}+\lambda v,
\end{equation}
for a some $\lambda$, or what is the same $Bv=\dot{v}+\lambda_jv$,
where $\lambda_j$ is the principal part of the Laurent series
expansion of $\lambda$ at $p$. Then given the curve $\mathcal{C}$
defined by (2.2) and $p\in\mathcal{C}$, Griffiths defines
\begin{eqnarray}
[\mbox{Laurent tail}(B)]_p&\equiv& \{\mbox{principal part of the Laurent series expansion}\nonumber\\
&&\mbox{ of }\lambda \mbox{ at }p\},\nonumber
\end{eqnarray}
and shows that the Lax flow can be linearized on the Jacobian
variety $\mbox{Jac}(\mathcal{C})$ if and only if for every
$p\in(h)_\infty$ (divisor of the poles of $h$), we have
\begin{eqnarray}
\frac{d}{dt}[\mbox{Laurent tail}(B)]_p&\in &\mbox{ linear combination } \{[\mbox{Laurent tail}(B)]_p \mbox{; Laurent}\nonumber\\
&&\mbox{ tail at } p \mbox{ of any meromorphic function }f \mbox{ on } \mathcal{C}\nonumber\\
&&\mbox{ such that } : (f)\geq n(h)_\infty\}.\nonumber
\end{eqnarray}
Equation (2.1) is invariant under the substitution $$B\longmapsto
B+P(h,A), \quad P(h,g)\in \mathbb{C}[h,g],$$ which shows that $B$
is not unique and that its natural place is somewhere in a
cohomology group. Let $$B(t,h)=\sum_{k=0}^nB_kh^k,$$ be a
polynomial of degree $n$. Let
$$\mathcal{D}=h^{-1}(\infty)=\sum_jn_jp_j, \quad n_j\geq 0,$$ (where $h$
is seen as a meromorphic function) be a positive divisor on
$\mathcal{C}$ and let $z_j$ a local coordinate around $p_j$. $B$
must be interpreted as an element of
$H^0(\mathcal{C},\mbox{Hom}(V,V(\mathcal{D}))$ where $V$ is the
sheaf of sections of the trivial bundle $\mathcal{C}\times V$ and
$V(\mathcal{D})=V\otimes \mathcal{O}_{\mathcal{C}}(\mathcal{D})$.
A section of $\mathcal{O}_{\mathcal{D}}(\mathcal{D})$ is written
$$\varphi=\sum\varphi_j, \quad
\varphi_j=\sum_{k=-n_j}^{-1}a_kz_j^k,$$ it is a principal part
(Laurent tail) centered on $p_j$. The Mittag-Leffler problem can
be formulated as follows : given a principal part $\varphi_j$,
find conditions for a function $\varphi\in
H^0(\mathcal{C},\mathcal{O}_{\mathcal{C}}(\mathcal{D}))$ such that
$\varphi-\varphi_j$ is holomorphic around $p_j$. The answer is
provided by the following theorem :

\begin{Theo}
Let $\mathcal{D}=\sum_ja_jp_j$. Given Laurent tail
$\{\varphi_j\}$, then there exist $\varphi\in
H^0(\mathcal{C},\mathcal{O}_{\mathcal{C}}(\mathcal{D}))$ such that
$\varphi-\varphi_j$ is holomorphic near $p_j$ if and only if
$$\sum_j\mbox{Res}_{p_j}(\varphi_j.\omega)=0,$$ for every
holomorphic differential $\omega$ on $\mathcal{C}$.
\end{Theo}

The residue of $B$, denoted by $\rho(B)\in
H^0(\mathcal{C},\mathcal{O}_{\mathcal{D}}(\mathcal{D})$, is the
collection of Laurent tails $\{\lambda_j\}$ given above (recall
that $\lambda_j$ is the principal part of the Laurent series
expansion of $\lambda$ at $p$). We shall say that the flow $L_t$
is linear if there exists a complex number $a$ such that
$$\frac{d^2L_t}{dt^2}=a\frac{dL_t}{dt}.$$ The Griffiths theorem is
as follows :

\begin{Theo}
We have
$$\dot{L}=\left.\frac{dL_t}{dt}\right|_{t=0}=\delta_1(\rho(B)).$$
Let $\mbox{Im res} \subset
H^0(\mathcal{C},\mathcal{O}_{\mathcal{D}}(\mathcal{D})$ be the
Laurent tails of meromorphic functions in
$H^0(\mathcal{C},\mathcal{O}_{\mathcal{D}}(\mathcal{D})$. Then the
flow $L_t$ (2.4) in $\mbox{Pic}^d(\mathcal{C})$ is linear if and
only if
\begin{equation}\label{eqn:euler}
\rho(\dot{B})=0 \mbox{ mod.} (\rho(B), \mbox{Im res}).
\end{equation}
\end{Theo}

The condition (2.8) is equivalent to
$$\sum_j\mbox{Res}_{p_j}(\dot{\rho}_j(B))\omega)=\mu\sum_j\mbox{Res}_{p_j}(\rho_j(B))\omega),\quad \omega\in
H^0(\mathcal{C},\Omega_{\mathcal{C}})$$ If this is satisfied, then
the linear flow on $\mbox{Jac}(\mathcal{C})$ is given by the
bilinear map
\begin{equation}\label{eqn:euler}
(t,\omega) \longmapsto
t\sum_j\mbox{Res}_{p_j}(\rho_j(B))\omega)=t\sum_j
\mbox{Res}_{p_j}(\lambda_j\omega).
\end{equation}

\begin{Exmp}
Let $J=\mbox{diag}(\lambda_1,...,\lambda_n)$, $\lambda_j>0$, be
the matrix representing the tensor of inertia of a $n$-dimensional
rigid body in a principal axis system and $\Omega(t)\in so(n)$ the
skew-symmetric matrix associated to the angular velocity vector of
the rigid body in the usual way. Define $M=\Omega J+J\Omega\in
so(n)$ the equations of motion of the rigid body can be written as
$$\dot{M}=[M,\Omega].$$ These equations are Hamiltonian on each
adjoint orbit of $so(n)$ defined by initial conditions with
Hamiltonian
$$H(M)=\frac{1}{2}(M,\Omega)=-\frac{1}{4}Tr(M\Omega).$$ By
Manakov's trick [43], these equations are equivalent to a Lax
equation with parameter
$$\dot{(\overbrace{M+J^2h})}=[M+J^2h, \Omega+Jh].$$
Hence $$\mathcal{D}=h^{-1}(\infty)=\sum_ip_i,$$ is the divisor
with $p_i$ being the $n$ distinct points lying over $h=\infty$. If
$z_i=h^{-1}$ is a local coordinate on $\mathcal{C}$ near $p_i$,
then from equation (2.7) with $B=\Omega+Jh$, we obtain
$$\rho(B)=\sum_i\frac{\lambda_i}{z_i}.$$ Since $\lambda_i$ are
constant, one has $\rho(\dot B)=0$ so the flow is linear on
$Jac(\mathcal{C}$). Since $$A=M+J^2h, \quad M+M^\top=0, \quad
J^2-{J^2}^\top=0,$$ we have $P(h,z)=(-1)^nP(-h,-z)$ and there is
an involution of the spectral curve
$$\sigma:\mathcal{C}\longrightarrow\mathcal{C},
\quad(h,z)\longmapsto(-h,-z).$$ We note that $\Omega$ moves on an
adjoint orbit $\mathcal{O}_\nu\subset so(n)$ and to linearize the
flow in question we need $\frac{1}{2}\dim\mathcal{O}_\nu$
integrals of motion that are in involution where for general
$\nu$,
\begin{equation}\label{eqn:euler}
\dim\mathcal{O}_\nu=\frac{n(n-1)}{2}-\left[\frac{n}{2}\right],
\end{equation}
Let $g(\mathcal{C})$ be the genus of the spectral curve
$\mathcal{C}$ and $g(\mathcal{C}_0)$ the genus of the quotient
$\mathcal{C}_0=\mathcal{C}/\sigma$ of $\mathcal{C}$ by the
involution $\sigma$. Since
$$g(\mathcal{C})=\frac{(n-1)(n-2)}{2},$$ then by  the
Riemann-Hurwitz formula,
\begin{equation}\label{eqn:euler}
g(\mathcal{C}_0)=g(\mathcal{C})-\frac{1}{2}\left(\frac{n(n-1)}{2}-\left[\frac{n}{2}\right]\right)=\left\{\begin{array}{rl}
\frac{(n-2)^2}{4}&n\equiv0\mbox{ mod.}2\\
\frac{(n-1)(n-3)}{4}&n\equiv1\mbox{ mod.}2
\end{array}\right.
\end{equation}
Associated to the double covering
$\mathcal{C}\longrightarrow\mathcal{C}_0$ is the Prym variety
$Prym(\mathcal{C}/\mathcal{C}_{0})$ and since
$\sigma(\rho(B))=-\rho(B)$, the flow in question actually occurs
on this complex torus. From (2.11) it follows that
\begin{equation}\label{eqn:euler}
\dim
Prym(\mathcal{C}/\mathcal{C}_{0})=\frac{1}{2}\left(\frac{n(n-1)}{2}-\left[\frac{n}{2}\right]\right)=\left\{\begin{array}{rl}
\frac{n(n-2)}{4}&n\equiv0\mbox{ mod.}2\\
\frac{(n-1)^2}{4}&n\equiv1\mbox{ mod.}2
\end{array}\right.
\end{equation}
On the other hand, comparing with (2.10) we obtain
$$\dim Prym(\mathcal{C}/\mathcal{C}_{0})=\frac{1}{2}\dim\mathcal{O}_\nu,$$
and we see that the motion of the free rigid linearizes on a torus
$Prym(\mathcal{C}/\mathcal{C}_{0})$ of exactly the right
dimension.
\end{Exmp}

\section{Continued fraction, orthogonal polynomials, the spectrum of Jacobi matrices and difference operators}

A Jacobi matrix is a doubly infinite matrix $(a_{ij})$ for $i
,j\in \mathbb{Z}$ such that : $a_{ij}=0$ if $|i-j|$ is large
enough. We show that the set of these matrices forms an
associative algebra and consequently a Lie algebra by
anti-symmetrization. Consider the Jacobi matrix
$$
\Gamma=\left(\begin{array}{ccccc}
b_{1}&a_{1}&0&\cdots&0\\
a_{1}&b_{2}&a_{2} &&\vdots \\
0&a_{2}&\ddots&\ddots&0\\
\vdots&&\ddots&\ddots&\\
0&\cdots&0&&
\end{array}\right),
$$
where all the $b_j$ are real and all the $a_j$ positive, and the
associated continued $\Gamma$-fraction,
\begin{equation}\label{eqn:euler}
\varphi(z)=\frac{a_0^2}{z-b_1-\displaystyle{\frac{a_1^2}{z-b_2-\displaystyle{\frac{a_2^2}{z-b_3-
_{\ddots}}}}}}
\end{equation}
where $a_0$ is a positive real number. By cutting off the
$\Gamma$-fraction $\varphi(z)$ at the $k$-th term, we obtain the
$k$-th Pad\'{e} approximant $\displaystyle{\frac{A_k(z)}{B_k(z)}}$
of $\varphi(z)$, i.e.,
\begin{equation}\label{eqn:euler}
\varphi(z)=\lim_{k\rightarrow\infty}\frac{A_k(z)}{B_k(z)}.
\end{equation}
The degree of the polynomial $A_k(z)$ is $k-1$, while the degree
of $B_(z)$ is $k$. Moreover, $\varphi(z)$ admits formal series
expansion in a neighborhood of the pole $z=0$, in the following
form
$$\varphi(z)=\frac{c_0}{z}+\frac{c_1}{z^2}+\frac{c_2}{z^3}+\cdots=\sum_{j=0}^\infty\frac{c_j}{z^{j+1}}.$$
Note that the characteristic polynomial $B_k$ of the
$\Gamma$-Jacobi matrix
$$B_k(z)=\det\left(\begin{array}{ccccc}
b_{1}-z&a_{1}&0&\cdots&0\\
a_{1}&b_{2}-z&a_{2} &&\vdots \\
0&a_{2}&\ddots&\ddots&0\\
\vdots&&\ddots&\ddots&a_{k-1}\\
0&\cdots&0&a_{k-1}&b_{k}-z
\end{array}\right),$$
is the last term of the second order recursion
$$B_j(z)-(z-b_k)B_{j-1}(z)+a_{j-1}^2B_{j-2}(z)=0.$$
The polynomials $A_k(z)$, $B_k(z)$ form a pair of linearly
independent solutions of a second order finite difference equation
(the eigenvectors of the Jacobi matrix from which we remove the
first row and column) :
$$a_{j}y_{j}+b_{j+1}y_{j+1}+a_{j+1}y_{j+2}=zy_{j+1},\quad j=0,1,...$$
with the boundary conditions : $$y_0\neq0, \quad y_1=0, \quad
y_{N+1}=0.$$ We have also the relation :
$$A_{j-1}(z)B_j(z)-A_j(z)B_{j-1}(z)=\frac{1}{a_{j-1}},\quad
j=1,2,...$$ From the classical theory, the polynomials $B_k$ form
an orthogonal system with respect to a Stieltjes measure
$d\sigma(x)$ on the real axis,
$$\int_{-\infty}^\infty B_k(x)B_l(x)d\sigma(x)=\delta_{kl}.$$
Conversely, if a family of polynomials $P_n(x)$ is orthogonal for
$d\sigma(x)$, then $P_n(x)$ satisfies the following recurrence
relation:
$$P_j(x)-(\lambda_j x-\mu_j)P_{j-1}(x)+\gamma_{j-1}P_{j-2}(x)=0,$$
where $\lambda_j>0$, $\mu$ and $\gamma_j>0$ are constants.
Moreover, if we consider the continued fraction
$$\psi(z)=\frac{\gamma_0}{\lambda_1z-\mu_1-\displaystyle{\frac{\gamma_1}{\lambda_2z-\mu_2-\displaystyle{\frac{\gamma_2}{\lambda_3z-\mu_3-
_{\ddots}}}}}}$$ and realize an equivalent transformation
$$\psi(z)=\frac{\gamma_0}{z-\frac{\mu_1}{\lambda_1}-\displaystyle{\frac{\frac{\gamma_1}{\lambda_1\lambda_2}}
{z-\frac{\mu_2}{\lambda_2}-\displaystyle{\frac{\frac{\gamma_2}{\lambda_2\lambda_3}}{z-\frac{\mu_3}{\lambda_3}-
_{\ddots}}}}}}$$ we reconstruct the $\Gamma$-fraction
corresponding to $d\sigma(x)$ (where we can put
$\frac{\gamma_j}{\lambda_j\lambda{j+1}}=a_j^2$ and
$\frac{\mu_j}{\lambda_j}=b_j$). It follows that there is a
one-to-one correspondence between the set of Jacobi matrices and
that of all the orthogonal polynomial systems on $\mathbb{R}$. In
fact, if we choose the orthogonal polynomials
$$P_n=\frac{\gamma_0}{\prod_{j=1}^{n-1}}B_{n-1}(x),$$ as the basis
of the vector space consisting of all polynomials then the Jacobi
matrix represents the multiplication by $x$.

As an example of $V_{-j,k}$ (theorem 3, b)), consider the infinite
Jacobi matrix (symmetric, tridiagonal and $N$-periodic) :
\begin{equation}\label{eqn:euler}
{A}=\left(\begin{array}{cccccccc}
\ddots&\ddots&&&&&&\\
\ddots&b_{0}&a_{0}&0&\cdots&0&\\
&a_{0}&b_{1}&a_{1} &&\vdots&\\
&0&a_{1}&\ddots&\ddots&0&&\\
&\vdots&&\ddots&\ddots&a_{N-1}&\\
&0&\cdots&0&a_{N-1}&b_{N}&\ddots\\
&&&&&\ddots&\ddots
\end{array}\right),\quad (a_i, b_i\in \mathbb{C})
\end{equation}
The matrix $A$ is $N$-periodic when $$a_{i+N}=a_i, \quad
b_{i+N}=b_i, \quad\forall i\in \mathbb{Z}.$$ We denote by
$f=(...,f_{-1},f_0,f_1,...)$ the (infinite) column vector and by
$D$ (shift operator) the operator passage of degree +1,
$Df_i=f_{i+1}$. Since the matrix $A$ is $N$-periodic, we have
$$AD^N=D^NA.$$ Reciprocally, this relation of commutation means that
$N$ is the period of $A$. Let
$$
{A}(h)=\left(\begin{array}{ccccc}
b_{1}&a_{1}&0&\cdots&a_{N}h^{-1}\\
a_{1}&b_{2}&a_{2} &&\vdots \\
0&a_{2}&\ddots&\ddots&0\\
\vdots&&\ddots&\ddots&a_{N-1}\\
a_{N}h&\cdots&0&a_{N-1}&b_{N}
\end{array}\right),\quad h\in \mathbb{C}^*
$$
be the finite Jacobi matrix (symmetric tridiagonal and
$N$-periodic). The determinant of the matrix
\begin{equation}\label{eqn:euler}
{A}(h)-zI=\left(\begin{array}{ccccc}
b_{1}-z&a_{1}&0&\cdots&a_{N}h^{-1}\\
a_{1}&b_{2}-z&a_{2} &&\vdots \\
0&a_{2}&\ddots&\ddots&0\\
\vdots&&\ddots&\ddots&a_{N-1}\\
a_{N}h&\cdots&0&a_{N-1}&b_{N}-z
\end{array}\right),
\end{equation}
is
\begin{equation}\label{eqn:euler}
F(h,h^{-1},z)\equiv\det
(A(h)-zI)=(-1)^{N+1}\left(\alpha\times(h+h^{-1})-P(z)\right),
\end{equation}
where $(z,h)\in \mathbb{C}\times \mathbb{C}^*$,
$\alpha=\prod_{i=1}^Na_i$ and $P(z)$ is a polynomial of degree $N$
with real coefficients :
\begin{eqnarray}
P(z)&=&\det\left(\begin{array}{ccccc}
b_{1}-z&a_{1}&0&\cdots&0\\
a_{1}&b_{2}-z&a_{2} &&\vdots \\
0&a_{2}&\ddots&\ddots&0\\
\vdots&&\ddots&\ddots&a_{N-1}\\
0&\cdots&0&a_{N-1}&b_{N}-z
\end{array}\right)\nonumber\\
&&\qquad-a_0^2\det\left(\begin{array}{ccccc}
b_{2}-z&a_{2}&0&\cdots&0\\
a_{2}&b_{3}-z&a_{3} &&\vdots \\
0&a_{3}&\ddots&\ddots&0\\
\vdots&&\ddots&\ddots&a_{N-2}\\
0&\cdots&0&a_{N-2}&b_{N-1}-z
\end{array}\right),\nonumber\\
&=&z^N+\cdots\nonumber
\end{eqnarray}
Let $\mathcal{C}$ be the Riemann surface defined by
\begin{eqnarray}
\mathcal{C}&=&\left\{(z,h)\in \mathbb{C}\times
\mathbb{C}^* : Af=zf, D^Nf=hf\right\},\nonumber\\
&=&\left\{(z,h)\in \mathbb{C}\times \mathbb{C}^* :
F(h,h^{-1},z)=0\right\}.
\end{eqnarray}
We suppose that $\alpha\neq 0$. From the equation $$F (h, h^{-1},
z)=0,$$ we derive the following relation : $$h=\frac{P(z)\pm
\sqrt{P^2(z)-4\alpha^2}}{2\alpha}.$$ Note that $\mathcal {C}$ is a
hyperelliptic curve branched in $2N$ points given by the roots of
the polynomial $P(z)=\pm2\alpha$, and admits two points at
infinity $\mathcal{P}$ and $\mathcal{Q}$; the point $\mathcal{P}$
covering the case $z=\infty$, $h=\infty$ while the point
$\mathcal{Q}$ is relative to the case $z=\infty$, $h=0$. (The
hyperelliptic involution on $\mathcal {C}$ maps $(z,h)$ into
$(z,h^{-1})$ and the curve $\mathcal {C}$ may be singular).
According to the Riemann-Hurwitz formula, the genus of
$\mathcal{C}$ is $g=N-1$. The meromorphic function $h$ has neither
zero nor poles except in the neighborhood of $z=\infty$. When
$z\nearrow\infty$, we have on the sheet +,
$$h\simeq\frac{P(z)+P(z)}{2\alpha}=\frac{P(z)}{\alpha}=\frac{z^N}{\alpha}+\cdots,$$
which shows that $h$ has a pole of order $N$. Similarly, when
$z\nearrow \infty$, we have on the sheet -,
$$h=\frac{P(z)-\sqrt{P^2(z)-4\alpha^2}}{2\alpha}=\frac{2\alpha}{P(z)+\sqrt{P^2(z)-4\alpha^2}}\simeq
\frac{\alpha}{z^N}+\cdots,$$ and therefore $h$ has a zero of order
$N$. Therefore the divisor $(h)$ of the function $h$ on
$\mathcal{C}$ is $$(h)=-N\mathcal{P}+N\mathcal{Q},$$ where
$\mathcal{P}$ and $\mathcal{Q}$ are the two points covering
$\infty$ on the sheets + and - respectively. The map $$\sim  :
\mathcal{C}\longrightarrow \mathcal{C}, \quad(z,h)\longmapsto
(\overline{z},\overline{h}^{-1}),$$ is an antiholomorphic
involution. In other words, the map $\sim :p\longmapsto
\widetilde{p}$ is such that :
$\widetilde{\mathcal{P}}=\mathcal{Q}$. When $|h|=1$, the finite
matrix $A(h)$ is self-adjoint and therefore admits a real
spectrum. Hence, the set of fixed points of this involution
denoted by $\mathcal{C}^\sim$ is determined by
$$
\mathcal{C}^\sim=\{p\in \mathcal{C}:
\widetilde{p}=p\}=\left\{(z,h): h=\overline{h}^{-1},
\overline{z}=z\right\}=\{(z,h): |h|=1\}.$$ Note that this set
divides $\mathcal{C}$ into two distinct regions $\mathcal{C}_ +$
and $\mathcal{C}_-$. More precisely, we have
$$
\mathcal{C}\backslash
\mathcal{C}^\sim=\mathcal{C}_+\cup\mathcal{C}_-=\{p\in \mathcal{C}
: |h|>1\}\cup \{p\in \mathcal{C} : |h|<1\},$$ so
$\mathcal{C}=\mathcal{C}_+\cup\mathcal{C}^\sim\cup\mathcal{C}_-$.
The first region $\mathcal{C}_+$ contains the point $\mathcal{P}$
while the second $\mathcal{C}_-$ contains the point $\mathcal{Q}$.
In fact $\mathcal{C}^\sim$ can be seen as the frontier of
$\mathcal{C}_+$ and $\mathcal{C}_-$, So $\mathcal{C}^\sim$ is
homologous to zero. Moreover, the involution $\sim$ extends to an
involution $*$ on the field of meromorphic functions as follows:
$$\varphi^*(p)=\overline{\varphi(\widetilde{p}}),$$ and on the
differential space as follows : $$(\varphi
d\psi)^*=\varphi^*d\psi^*.$$ Hence, we have $h^*=h^{-1}$ and
$z^*=z$. The condition that the matrices $A$ and $D^N$ have an
eigenvector in common is parameterized by the Riemann surface
$\mathcal{C}$ (2.6), let $f=(...,f_{-1},f_0,f_1,...)$ such an
eigenvector. In the following, appropriate standardization is used
by selecting $f_0\equiv 1$, from where $F_N=h$. Let us therefore
$\overline{f}=(f_1,f_2,...,f_{N-1})^\top$. Since $\overline{f}$
satisfies  $$(A (h)-zI)\overline{f}=0,$$ then we have
$$f_k=\frac{\Delta_{1,k}}{\Delta_{1,l}}f_l=\frac{\Delta_{2,k}}{\Delta_{2,l}}f_l=\cdots
=\frac{\Delta_{N,k}}{\Delta_{N,l}}f_l,\quad 1\leq k,l\leq N,$$
where $\Delta_{k,l}$ is the $(k,l)$-cofactor of $(A(h)-zI)$, that
is to say,
\begin{equation}\label{eqn:euler}
\Delta_{k,l}=(-1)^{k+l}\times (k,l)-\mbox{minor of }
(A(h)-zI).
\end{equation}
is the $(k,l)-\mbox{minor of }(A(h)-zI)$, i.e., the determinant of
the $N-1$ submatrix obtained by removing the $k^{th}$-line and the
$l^{th}$-column of the matrix $(A(h)-zI)$). In particular, $f$ can
be expressed as a rational function in $z$ and $h$,
$$f_k=\frac{\Delta_{N,k}}{\Delta_{N,N}}h=\frac{\Delta_{k,k}}{\Delta_{k,N}}h.$$
According to matrix (3.4), we note that
\begin{eqnarray}
\Delta_{N,1}&=&\prod_{j=1}^{N-1}a_j+(-1)^Na_Nh^{-1}\left((-z)^{N-2}+\cdots\right),\nonumber\\
\Delta_{1,N}&=&\prod_{j=1}^{N-1}a_j+(-1)^Na_Nh\left((-z)^{N-2}+\cdots\right),\nonumber
\end{eqnarray}
where $(-z)^{N-2}+\cdots,$ is a polynomial of degree $N-2$.
Similarly,
$$
\Delta_{N,N}=(-z)^{N-1}+\cdots,$$ is a polynomial of degree $N-1$.
To determine the divisor structure of $f_k$, one proceeds as
follows : for $f_1$, we have
\begin{eqnarray}
(f_1)_\infty&=&(\Delta_{N,1})_\infty+(h)-(\Delta_{N,1})_\infty,\nonumber\\
&=&-(2N-2)\mathcal{Q}-N\mathcal{P}+N\mathcal{Q}+(N-1)\mathcal{P}+(N-1)\mathcal{Q},\nonumber\\
&=&\mathcal{Q}-\mathcal{P},\nonumber
\end{eqnarray}
and for the other $f_k$, we consider first the matrix (3.4)
shifted by one, i.e.,
$$
\left(\begin{array}{ccccc}
b_{2}-z&a_{2}&0&\cdots&a_{1}h^{-1}\\
a_{2}&b_{3}-z&a_{3} &&\vdots \\
0&a_{3}&\ddots&\ddots&0\\
\vdots&&\ddots&b_{N}-z&a_{N}\\
a_{1}h&\cdots&0&a_{N}&b_{1}-z
\end{array}\right).
$$
Hence,
$$
\left(\begin{array}{ccccc}
b_{2}-z&a_{2}&0&\cdots&a_{1}h^{-1}\\
a_{2}&b_{3}-z&a_{3} &&\vdots \\
0&a_{3}&\ddots&\ddots&0\\
\vdots&&\ddots&b_{N}-z&a_{N}\\
a_{1}h&\cdots&0&a_{N}&b_{1}-z
\end{array}\right)\left(\begin{array}{c}
\frac{f_2}{f_1}\\
\frac{f_3}{f_1}\\
\vdots\\
\frac{f_N}{f_1}\\
h
\end{array}\right)=0,
$$
and as above, we have
$$\left(\frac{f_2}{f_1}\right)_\infty=\mathcal{Q}-\mathcal{P},$$
which implies that
$$(f_2)_\infty=\left(\frac{f_2}{f_1}\right)_\infty+(f_1)_\infty=2\mathcal{Q}-2\mathcal{P}.$$
And in general, we get
$$(f_k)_\infty=k\mathcal{Q}-k\mathcal{P}.$$ Let $\mathcal {D}$ be
a minimal positive divisor on $\mathcal{C}$ such that :
$$(f_k)+\mathcal{D}\geq-k\mathcal{P}+k\mathcal{Q}, \quad\forall k\in
\mathbb{Z}.$$ It is shown that the degree of $\mathcal{D}$ is
$$\mbox{deg }\mathcal{D}=g=N-1.$$ We show that the divisor
$\mathcal{D}$ is regular with respect to $\mathcal{P}$ and
$\mathcal{Q}$, i.e., such that $$\dim
\mathcal{L}(\mathcal{D}+k\mathcal{P}-(k+1)\mathcal{Q})=0,\quad
\forall k\in \mathbb{Z}$$ The proof consists in showing first that
the divisor $\mathcal{D}$ is general. Recall that a positive
divisor $\mathcal{D}$ of degree $g$ on $\mathcal{C}$ is general if
$(\omega_j(p_j))\neq 0$, $p_j\in \mathcal{C}$, $1\leq j\leq g$
where $(\omega_1,...,\omega_g)$ is a normalized base of
differential forms on $\mathcal{C}$. It is shown that
$\mathcal{C}$ is general if and only if
$\dim\mathcal{L}(\mathcal{D})=1$ (where $\mathcal{L}(\mathcal{D})$
denotes the set of meromorphic functions $f$ on $\mathcal{C}$ such
that : $(f)+\mathcal{D}\geq0$) or if and only if
$\dim\Omega(-\mathcal{D})=0$ where $\Omega(\mathcal{D})$ denotes
the set of meromorphic differential forms $\omega$ on
$\mathcal{C}$ such that the divisor $(\omega)+\mathcal{D}\geq0$.
Consider an integer $k>g-2$, then we deduce from the Riemann-Roch
theorem
$$\dim\mathcal{L}(\mathcal{D}+k\mathcal{P})=\dim\Omega(-\mathcal{D}-k\mathcal{P})+g+k-g+1=\dim\Omega(-\mathcal{D}-k\mathcal{P})+k+1.$$
Since $\dim\Omega(-\mathcal{D}-k\mathcal{P})=0$, because a
holomorphic differential can have at most $2g-2$ zeroes, then
$$\dim\mathcal{L}(\mathcal{D}+k\mathcal{P})=k+1.$$ Moreover,
$\mathcal{L}(\mathcal{D}+j\mathcal{P})$ is strictly larger than
$\mathcal{L}(\mathcal{D}+(j-1)\mathcal{P})$, because $f_j$ belongs
to the first space and not to the second. Therefore by lowering
the index $j$ down to $0$, it follows that
$\dim\mathcal{L}(\mathcal{D})=1$, which shows that $\mathcal{D}$
is general. Let's show now that $\mathcal{D}$ is regular. It
suffices to proceed by induction. Since
$\dim\mathcal{L}(\mathcal{D})=1$ and
$\mathcal{L}(\mathcal{D}-\mathcal{Q})\subsetneqq\mathcal{L}(\mathcal{D})$
(i.e., $f_0=1\notin\mathcal{L}(\mathcal{D}-\mathcal{Q})$; the
function $f_0=1$ belongs to the second space but not the first),
then $$\dim\mathcal{L}(\mathcal{D}-\mathcal{Q})=0.$$ Assume that
$$\dim
\mathcal{L}(\mathcal{D}+k\mathcal{P}-(k+1)\mathcal{Q})=0,$$ then
by the Riemann-Roch theorem
$$\dim \mathcal{L}(\mathcal{D}+(k+1)\mathcal{P}-(k+2)\mathcal{Q})\leq
\dim \mathcal{L}(\mathcal{D}+k\mathcal{P}-(k+1)\mathcal{Q})+1=1,$$
implies equality since $f_{k+1}$ belongs to the first space. Since
$f_{k+1}$ belongs to
$\mathcal{L}(\mathcal{D}+(k+1)\mathcal{P}-(k+1)\mathcal{Q})$ but
not to
$\mathcal{L}(\mathcal{D}+(k+1)\mathcal{P}-(k+2)\mathcal{Q})$, we
have that
$$\dim\mathcal{L}(\mathcal{D}+(k+1)\mathcal{P}-(k+2)\mathcal{Q})=0.$$
Consider the differential of $F$ (2.34) while taking into account
that $z$ appears only on the diagonal of the matrix $A(h)-zI$. We
have $$-\sum_{i=1}^N\Delta_{ii}dz+h\frac{\partial F}{\partial
h}\frac{dh}{h}=0,$$ and either
$$\omega=\frac{-i\Delta_{NN}dz}{h\frac{\partial F}{\partial
h}}.$$ We have
$$
\omega=\frac{-i\frac{dh}{h}}{\sum_{i=1}^N\frac{\Delta_{ii}}{\Delta_{NN}}}=
\frac{-i\frac{dh}{h}}{\sum_{i=1}^N\frac{\Delta_{ii}}{\Delta_{iN}}.\frac{\Delta_{iN}}{\Delta_{NN}}}
=\frac{-i\frac{dh}{h}}{\sum_{i=1}^N\frac{\Delta_{Ni}}{\Delta_{NN}}.\frac{\Delta_{iN}}{\Delta_{NN}}}.$$
Or $$\Delta_{iN}=\Delta_{Ni}^*, \quad1\leq i\leq N,$$ so
$$
\omega=\frac{-i\frac{dh}{h}}{\sum_{i=1}^N\frac{\Delta_{Ni}}{\Delta_{NN}}
.\left(\frac{\Delta_{iN}}{\Delta_{NN}}\right)^*}=\frac{-i\frac{dh}{h}}{\sum_{i=1}^Nf_if_i^*},$$
and consequently
$$\omega=\pm\frac{\Delta_{NN}dz}{\sqrt{P^2(z)-4A^2}}.$$
From this we deduce that $\omega^*=\omega$. In addition,
$\omega\geq 0$ on $\mathcal{C}^\sim$. (Indeed, on
$\mathcal{C}^\sim$ we have $$\sum_{i=1}^Nf_if_i^*
=\sum_{i=1}^N|f_i|^2\geq 0.$$ Let $h=\rho e^{i\theta}$. Note that
in all finite number points, $h$ is a local parameter on
$\mathcal{C}$ while $\theta$ is a local parameter on
$\mathcal{C}^\sim$. Like $$-ih^{-1}dh=d\theta, \quad\omega\geq
0,$$ at these points and by continuity at all points). We also
have a relation which shows that the scalar product between $f_k$
and $f_l$ is
$$\langle f_k,
f_l\rangle=\int_{\mathcal{C}^\sim}f_k.f_l^*\omega=\left\{\begin{array}{rl}
0&\mbox{si }k\neq l\\
>0&\mbox{si }k=l
\end{array}\right.
$$
That is, the functions $f_k$, $k\in \mathbb{Z}$, are orthogonal to
$\mathcal{C}^\sim$ with respect to $\omega$. We deduce from these
properties that the divisor of $\omega$ is
$$(\omega)=\mathcal{D}+\widetilde{\mathcal{D}}-\mathcal{P}-\mathcal{Q},$$
for the involution $\sim$ introduced previously. Given a matrix of
the form $A$ (2.3), we have obtained a series of data
$\{\mathcal{C}, z, h, \mathcal{D}, \omega\}$.

What is remarkable is that the reverse is also true (for further
information, see [60]) :

\begin{Theo}
There is a one-to-one correspondence between the following sets of
data :

a) Let $a_i, b_i\in \mathbb{C}$, $a_i\neq0$, where
$$a_{i+N}=a_i, \quad b_{i+N}=b_i, \quad -\infty<i<+\infty.$$ An infinite
$N$-periodic matrix
$$
\left(\begin{array}{cccccccc}
\ddots&\ddots&&&&&&\\
\ddots&b_{0}&a_{0}&0&\cdots&0&\\
&\overline{a}_{0}&b_{1}&a_{1} &&\vdots&\\
&0&\overline{a}_{1}&\ddots&\ddots&0&&\\
&\vdots&&\ddots&\ddots&a_{N-1}&\\
&0&\cdots&0&\overline{a}_{N-1}&b_{N}&\ddots\\
&&&&&\ddots&\ddots
\end{array}\right),
$$
modulo conjugation by $N$-periodic diagonal matrices with real
entries.

b) A curve $\mathcal{C}$ (possibly singular) of genus $N-1$ with
two points $\mathcal{P}$ and $\mathcal{Q}$ on $\mathcal{C}$, a
divisor $\mathcal{D}$ of degree $N-1$ on $\mathcal{C}$ and two
meromorphic functions $h$ and $z$ on $\mathcal{C}$ such that :
$$(h)=-N\mathcal{P}+N\mathcal{Q},\qquad (z)=-\mathcal{P}-\mathcal{Q}+S,$$
where $S$ is a positive divisor not containing the points
$\mathcal{P}$ and $\mathcal{Q}$. The curve $\mathcal{C}$ is
equipped with an antiholomorphic involution $$\sim
:(z,h)\longmapsto (\overline{z},\overline{h}^{-1}),$$ for which
$$\mathcal{C}=\mathcal{C}_+\cup\mathcal{C}^\sim\cup\mathcal{C}_-,$$
where $$\mathcal{C}^\sim=\{p\in
\mathcal{C}:\widetilde{p}=p\}=\{(z,h):|h|=1\},$$
$$\mathcal{C}_+=\{p \in \mathcal{C}:|h|>1\},$$ $$\mathcal{C}_-=\{p
\in \mathcal{C}:|h|<1\},$$ such that : $\mathcal{P}\in
\mathcal{C}_+$ and $\mathcal{Q}\in \mathcal{C}_-$. By using the
involution $\sim$, we introduce an involution $*$ acting on the
space of all functions on $\mathcal{C}$ in a way
$$\varphi^*(p)=\widetilde{\varphi(\widetilde{p}}),$$ and on the
differential space as follows : $$(\varphi
d\psi)^*=\varphi^*d\psi^*,$$ then $h^*=h^{-1}$, and $z^*=z$. The
divisor of a differential form $\omega$ on $\mathcal{C}$ is
$$(\omega)=\mathcal{D}+\widetilde{\mathcal{D}}-\mathcal{P}-\mathcal{Q}.$$
\end{Theo}

For any difference operator $X$, we define
$$
\left(X^{[+]}\right)_{ij}=\left\{\begin{array}{rl} X_{ij}&\mbox{si
}i<j,\\
\frac{1}{2}X_{ij}&\mbox{si }i=j,\\
0&\mbox{si }i>j,
\end{array}\right.,\qquad
X^{[-]}=X-X^{[+]}.
$$
Let $\mathcal{M}$ be the vector space of infinite $N$-periodic
matrices $A$ such that for some $K$, $c_{ij}=0$ if $|i-j|>K$. On
$\mathcal{M}$, we introduce the following scalar product :
$$\langle C,D\rangle=Tr(CD^\top)=\sum_{(i,j)\in
\mathbb{Z}^2}c_{ij}d_{ij}.$$ We say that a functional $F$ is
differentiable if there exists a matrix $\frac{\partial
F}{\partial C}$ in $\mathcal{M}$ such that for every $D$,
$$\lim_{\epsilon\rightarrow 0}\frac{F(C+\epsilon
D)-F(C)}{\epsilon}=\left\langle \frac{\partial F}{\partial
C},D\right\rangle.$$ Note that we also have
$$\langle[A,B],C\rangle=\langle[A^\top,C],B\rangle.$$ Define the
following bracket between two differentiable functionals $F$ and
$G$ on $\mathcal{M}$,
$$\{F,G\}=\left\langle \left[\left(\frac{\partial
F}{\partial X}\right)^{[+]},\left(\frac{\partial G}{\partial
X}\right)^{[+]} \right]-\left[\left(\frac{\partial F}{\partial
X}\right)^{[-]},\left(\frac{\partial G}{\partial X}\right)^{[-]}
\right], X\right\rangle.$$ $\{,\}$ satisfies the Jacobi identity.
Let $P(A, S, S^{-1})$ be a polynomial in $S+S^{-1}$ and $A$ with
real coefficients. Consider the following Lax equation:
\begin{equation}\label{eqn:euler}
\dot{A}=\left[P(A,S,S^{-1})^{[+]}-P(A,S,S^{-1})^{[-]},A\right].
\end{equation}
When the matrix $A(t)$ deforms with $t$, then only the divisor
$\mathcal{D}$ varies while
$\{\mathcal{C},z,h,\mathcal{P},\mathcal{Q}\}$ remain fixed. As we
have already shown, the coefficients of $z^ih^j$ in equation (2.5)
are invariants of this motion. The divisor $\mathcal{D}(t)$
evolves linearly on the Jacobian manifold
$\mbox{Jac}(\mathcal{C})$. Any linear flow over
$\mbox{Jac}(\mathcal{C})$ is equivalent to Equation (2.8) and is a
Hamiltonian flow with respect to the above (Poisson) bracket. In
particular, the flow
$$\dot{A}=\left[A,(S^{-k}A^l]^{[+]}\right],$$ is written in terms
of the (Poisson) bracket as follows :
$$\dot{a}_{ij}=\{F,a_{ij}\},\qquad F=\frac{1}{l+1}Tr\left(S^{-k}A^{l+1}\right).$$
The (Poisson) bracket of two functional of the form
$Tr\left(S^{-k}A^{l+1}\right)$ is zero, which means that we have a
set of integrals in involution. Let $(\omega_1, ...,\omega_g)$ be
a holomorphic differential basis on the hyperelliptic curve
$\mathcal{C}$. We have
$$\omega_k=\frac{z^{k-1}}{\sqrt{P^2(z)-4Q^2}},$$ and let
$$c_k=\mbox{Res}_p(\omega_kz^j), \quad1\leq j\leq g.$$ Since the order
of the zeros of $\omega_k$ at the points at infinity
$\mathcal{P}$, $\mathcal{Q}$ is equal to $g-k$, then $c_k=0$ for
$k<g-j+1$ and $c_k\neq 0$ for $k=g-j+1$. Therefore, a complete set
of flows is given by the functions $z,z^2,...,z^g$ and the flow
which leaves invariant the spectrum of $A$ and $X$ is given by a
polynomial $P(z)$ of degree at most equal to $g$ :
$$\dot{A}=\frac{1}{2}\left[A, P(A)^+-P(A)^-\right],$$ where $P(A)^+$
(resp. $-P(A)^-$) is the upper triangular part of $P(A)$ (lower),
including the diagonal of $P(A)$. The (Poisson) bracket between
two functional $F$ and $G$ can still be written in the form
$$\{F, G\}=\left\langle \left(\begin{array}{c}
\displaystyle{\frac{\partial F}{\partial a}}\\
\displaystyle{\frac{\partial F}{\partial b}}
\end{array}\right)^\top, J\left(\begin{array}{c}
\displaystyle{\frac{\partial G}{\partial a}}\\
\displaystyle{\frac{\partial G}{\partial b}}
\end{array}\right)\right\rangle,$$
where $\frac{\partial F}{\partial a}$ and $\frac{\partial
F}{\partial b}$ are the column vectors whose elements are given by
$\frac{\partial F}{\partial a_i}$ and $\frac{\partial F}{\partial
b}_i$ respectively while $J$ is the antisymmetric matrix Of order
$2n$ defined by
$$
{J}=\left(\begin{array}{cc}
O&\mathcal{A}\\
-\mathcal{A}^\top&O
\end{array}\right),
\qquad {\mathcal{A}}=2\left(\begin{array}{ccccc}
a_1&0&0&...&-a_N\\
-a_1&a_2&0&&\vdots\\
0&-a_2&a_3&&\vdots\\
\vdots&&&&\vdots\\
0&...&...&-a_{N-1}&a_N
\end{array}\right).
$$
The symplectic structure is given by
\begin{equation}\label{eqn:euler}
\omega=\sum_{j=2}^Ndb_j\wedge \sum_{j\leq i\leq
N}\frac{da_i}{a_i}.
\end{equation}
Flaschka variables [13] : $$a_{j}=\frac{1}{2}e^{x_{j}-x_{j+1}},
\qquad b_{j}=-\frac{1}{2}y_{j},$$ applied to the form (3.9) with
$x_{N+1}=0$, leads to the symplectic structure
$$\omega=\frac{1}{2}\sum_{j=2}^Ndx_j\wedge dy_j,$$
used by Moser [45, 46] during the study of the dynamical system of
$N-1$ particules moving freely on the real axis under the
influence of the exponential potential. See also the example below
concerning the study of Toda lattice. We have
$$\det (A_h-zI)|_{h=i}=(-1)^Nz^N+\sum_{i=1}^N\beta_iz^{N-1},$$
where $\beta_2,...,\beta_N$ are the $g$ invariant, functionally
independent and in involution. These invariants are given by the
$g=N-1$ points chosen from the spectrum of $A_1$ and $A_{-1}$,
i.e., by the branch points of the hyperelliptic curve
$\mathcal{C}$ or by the quantities $Tr A^k$, $2\leq k \leq N$.

With the Jacobi matrix, we can associate an operator $T$ on a
separable Hilbert space $E$ as follows,
\begin{eqnarray}
Te_0&=&b_0e_0+a_0e_1,\nonumber\\
Te_j&=&a_{j-1}e_{j-1}+b_je_j+a_je_{j+1},\quad j=1,2,... \nonumber
\end{eqnarray}
where $(e_1,...)$ is an orthonormal basis in $E$. The operator $T$
is symmetric. Indeed, for any two finite vectors $u$ and $v$, we
have $\langle Tu,v\rangle=\langle u,Tv\rangle$, according to the
symmetry of the Jacobi matrix. Moreover, if the Carleman's
condition : $$\sum_{j=0}^\infty\frac{1}{a_j}=+\infty$$ is
satisfied, then the operator $T$ is self-adjoint and its spectrum
is simple with $e_0$ a generating element. In this case, the
information about the spectrum of $T$ is contained in function,
\begin{equation}\label{eqn:euler}
\varphi(z)=\left\langle
(T-zI)^{-1}e_0,e_0\right\rangle=\int_{-\infty}^\infty\frac{d\sigma(x)}{z-x},
\end{equation}
defined at $z\notin\sigma(T)$ where $\sigma(x)=\langle
I_xe_0,e_0\rangle$ and $I_x$ is the resolution of the identity of
the operator $T$. Recall that the infinite continued fraction
converges if the limit (2.2) exists. If the operator $T$ is
self-adjoint, then the continued fraction $\varphi(z)$ converges
uniformly in any closed bounded domain of $z$ without common
points with real axis, to the analytic function defined by (3.10).
If the support of $d\sigma(x)$ is bounded, then the sequence
$\left(\frac{A_k(z)}{B_k(z)}\right)$ converges uniformly to a
holomorphic function near $z=\infty$. Moreover, if a Jacobi matrix
is bounded, i.e., if there exists $\rho>0$ such that, $$\forall j,
\quad |a_j|\leq\frac{\rho}{3}, \quad |b_j|\leq\frac{\rho}{3},$$
then the associated $\Gamma$-fraction converges uniformly on the
following domain $\{z : |z|\geq\rho\}$ and the support of
$d\sigma(x)$ is included in $[-\rho,\rho]$. In the case of a
periodic Jacobi matrix, this one is obviously bounded and
therefore the associated $\Gamma$-fraction converges near
$z=\infty$. In addition, the function $\varphi(z)$ is written in
the form (3.10) (Cauchy-Stieltjes transform of $d\sigma(x)$),
which shows that $\varphi(z)$ has zero of first order at
$z=\infty$ and for any point $z$ belonging the upper-half plane,
the imaginary part of $\varphi(z)$ is non positive.

We will now extend the Jacobi matrix $\Gamma$ to  the infinite
symmetric, tridiagonal and $N$-periodic Jacobi matrix $A$(3.3) and
use the results obtained previously. We consider $\varphi(z)$(3.1)
as being the associated $N$-periodic $\Gamma$-fraction. The latter
converges near the infinite point $z=\infty$. After an analytic
prolongation, the function $\varphi(z)$ coincides with $a_0f_1$
where $f_1$ is a meromorphic function on the genus $N-1$
hyperelliptic curve $\mathcal{C}$(3.6). The latter is branched at
the $2N$ real zeroes $\xi_1$, $\xi_2$,...,$\xi_{2N}$ of the
polynomial $P^2(z)-4\alpha^2$. The interval
$[\xi_{2j-1},\xi_{2j}]$, $1\leq j\leq N$, is called the stable
band and the interval $[\xi_{2j},\xi_{2j+1}]$, $1\leq j\leq N-1$,
is called the unstable band.

\begin{Theo}
Each zero $\sigma_1<\sigma_2<\cdots<\sigma_{N-1}$ of
$\Delta_{k,l}$ (3.7), belongs to the $j$-th finite unstable band
$[\lambda_{2j},\lambda_{2j+1}]$, $1\leq j\leq N-1$.
\end{Theo}

The function $\varphi(z)$ can be expressed (see below) by means of
Abelian integrals on the hyperelliptic curve $\mathcal{C}$ (3.6).
For $N=1$, $B_k(x)$ is the $k$-th Tschebyscheff polynomial of the
second type. For $N>1$, Kato [25, 26] has found a new phenomenon
related to discrete measures. We have seen that
$$\varphi(z)=a_0f_1=a_0\frac{\Delta_{N,1}}{\Delta_{N,N}}h,$$
belonging to $\mathcal{L}(\mathcal{D}+\mathcal{P}-\mathcal{Q})$.
Then, we have

\begin{Theo}
The function $\varphi(z)$ can explicitly written by means of
Abelian integrals on the hyperelliptic curve $\mathcal{C}$ (3.6)
as follows,
\begin{equation}\label{eqn:euler}
\varphi(z)=\sum_{j=1}^{N-1}\frac{\mbox{Res}(\varphi(z),\sigma_j^-)}{z-\sigma_j}+\sum_{j=1}^N\frac{(-1)^{N+1}}{2\pi
i}\int_{\xi_{2l-1}}^{\xi_{2l}}\frac{\sqrt{P^2(x)-4\alpha^2}}{(z-x)\Delta_{N,N}(x)}dx,
\end{equation}
where,
$$\mbox{Res}(\varphi(z),\sigma_j^-)=\frac{\alpha h(\sigma_j^-)+(-1)^Na_0^2.\Lambda}{\prod_{l\neq j}(\sigma_j-\sigma_l)},$$
and
$$\Lambda\equiv\det\left(\begin{array}{ccccc}
b_{2}-\sigma_j&a_{2}&0&\cdots&0\\
a_{2}&b_{3}-\sigma_j&a_{3} &&\vdots \\
0&a_{3}&\ddots&\ddots&0\\
\vdots&&\ddots&\ddots&a_{N-2}\\
0&\cdots&0&a_{N-2}&b_{N-1}-\sigma_j
\end{array}\right).$$
\end{Theo}

The differentials obtained in the previous section,
$$a\frac{\Delta_{N,N}(x)}{\sqrt{P^2(x)-4\alpha^2}}dx,\qquad
b\frac{\sqrt{P^2(x)-4\alpha^2}}{\Delta_{N,N}(x)}dx,$$ ($a$ and $b$
are constants) are positive mesures on each stable band
$[\xi_{2j-1},\xi_{2j}]$. Therefore, the expression (3.11) means
that $\varphi(z)$ can be obtained by the Cauchy-Stieltjes
transform of
\begin{eqnarray}
d\sigma&=&\sum_{j=1}^{N-1}\mbox{Res}(\varphi(z),\sigma_j^-).\delta(x-\sigma_j)dx+\frac{(-1)^{N+1}}{2\pi
i}.\frac{\sqrt{P^2(x)-4\alpha^2}}{\Delta_{N,N}(x)}dx,\nonumber\\
&=&\mbox{discrete mesure }+\mbox{ continuous mesure},\nonumber
\end{eqnarray}
as follows,
$$\varphi(z)=\int_{-\infty}^\infty\frac{d\sigma}{z-x}.$$
The function $\varphi(z)$ belongs to
$\mathcal{L}(\mathcal{D}'+\mathcal{P}-\mathcal{Q})$ where
$\mathcal{D}'=\sigma_1^++\cdots+\sigma_{N-1}^+$ is contained in
$\mathcal{C}_+=\{p\in \mathcal{C} : |h|>1\}$ (see previous
section). From expression (3.11), we have
$$\mathcal{D}=\sigma_{j_1}^-+\cdots+\sigma_{j_l}^-+\sigma_{j_{l+1}}^++\cdots+\sigma_{j_{N-1}}^+,$$
where $j_1<j_2<...<j_l$ denote the numbers for which
$\mbox{Res}(\varphi(z),\sigma_j^-)>0$ and
$j_{l+1}<j_{l+2}<...<j_{N-1}$ the numbers for which
$\mbox{Res}(\varphi(z),\sigma_j^-)=0$.
Hence,$$\mbox{Res}(\varphi(z),\sigma_j^-)=0 \mbox{ or }
-\frac{\sqrt{P^2(\sigma_{j}^-)-4\alpha^2}}{\prod_{l\neq
j}(\sigma_j-\sigma_l)}.$$

\begin{Exmp}
The Toda lattice equations [57] (discretized version of the
Korteweg-de Vries equation) describe the motion of $n$ particles
with exponential restoring forces and are governed by the
following Hamiltonian
$$H=\frac{1}{2}\sum_{j=1}^{N}y_{j}^{2}+\sum_{j=1}^{N}e^{x_{j}-x_{j+1}}.$$
Here the nonlinear dynamical system is described by the following
Hamiltonian equations :
$$
\dot{x}_{j}=y_{j},\qquad
\dot{y}_{j}=-e^{x_{j}-x_{j+1}}+e^{x_{j-1}-x_{j}}.$$ Flaschka
variables [13] :
$$a_{j}=\frac{1}{2}e^{x_{j}-x_{j+1}},\qquad b_{j}=-\frac{1}{2}y_{j},$$
can be used to express the symplectic structure $\omega$ (3.9) in
terms of $x_j$ and $y_j$ as follows,
$$
\frac{d{a}_{j}}{a_j}=dx_{j}-dx_{j+1},\qquad 2d{b}_{j}=-dy_{j},$$
then
$$
\omega=-\frac{1}{2}\sum_{j=2}^Ndy_j\sum_{i=j}^N(dx_i-dx_{i+1})=\frac{1}{2}\sum_{j=2}^Ndx_j^*\wedge
dy_j^*,
$$
where $x_j^*\equiv x_j-x_1$ and $y_j^*\equiv y_j$. We will study
the integrability of this problem with the Griffiths approach.
There are two cases :

$(i)$ The non-periodic case, i.e., $$x_{0}=-\infty, \quad
x_{N+1}=+\infty,$$ where the masses are arranged on a line. In
term of the Flaschka variables above, Toda's equations take the
following form
$$
\dot{a}_{j}=a_{j}\left( b_{j+1}-b_{j}\right),\qquad
\dot{b}_{j}=2(a_{j}^2-a_{j+1}^2),$$ with $a_{N+1}=a_{1}$ and
$b_{N+1}=b_{1}$. To show that this system is completely
integrable, one should find $N$ first integrals independent and in
involution each other. From the second equation, we have
$$\left(\sum_{j=1}^{N}b_{j}\right)^.=\sum_{j=1}^{N}\dot{b}_{j}=0,$$ and we normalize the $b_i$'s by requiring that
$\sum_{j=1}^{N}b_{j}=0$ (applying this fact to (3.9), leads to the
original symplectic form $\omega=\frac{1}{2}\sum_{j=2}^Ndx_j\wedge
dy_j$). This is a first integral for the system and to show that
it is completely integrable, we must find $N-1$ other integrals
that are functionally independent and in involution. We further
define $N\times N$ matrices $A$ and $B$ with
$$
{A}=\left(\begin{array}{ccccc}
b_{1}&a_{1}&0&\cdots&a_{N}\\
a_{1}&b_{2}&\vdots &&\vdots \\
0&\ddots&\ddots&\ddots&0\\
\vdots&&\ddots&b_{N-1}&a_{N-1}\\
a_{N}&\cdots&0&a_{N-1}&b_{N}
\end{array}\right),$$
$${B}=\left(\begin{array}{ccccc}
0&a_{1}&\cdots&\cdots&-a_{N}\\
-a_{1}&0&\vdots &&\vdots \\
\vdots&\ddots&\ddots&\ddots&\vdots\\
\vdots&&\ddots&\ddots&a_{N-1}\\
a_{N}&\cdots&\cdots&-a_{N-1}&0
\end{array}\right).
$$
Then the proposed system is equivalent to the Lax equation $\dot
A=[B,A]$. From theorem 1, we know that the quantities
$$I_k=\frac{1}{k}trA^{k},\quad1\leq k\leq N,$$ are first integrals of
motion. To be more precise
$$\dot I_k=tr(\dot A.A^{k-1})=tr([B,A].A^{k-1})
=tr(BA^k-ABA^{k-1})=0.$$ Notice that $I_1$ is the first integral
already know. Since these $N$ first integrals are shown to be
independent and in involution each other, the system in question
is thus completely integrable.

$(ii)$ The periodic case, i.e., $$y_{j+N}=y_{j}, \quad
x_{j+N}=x_{j},$$ the connected masses will be arranged on a
circle. We show that in this case, the spectrum of the periodic
Jacobi matrix
$$
A(h)=\left(\begin{array}{ccccc}
b_{1}&a_{1}&0&\cdots&a_{N}h^{-1}\\
a_{1}&b_{2}&\vdots &&\vdots \\
0&\ddots&\ddots&\ddots&0\\
\vdots&&\ddots&b_{N-1}&a_{N-1}\\
a_{N}h&\cdots&0&a_{N-1}&b_{N}
\end{array}\right),
$$
remains invariant in time. The matrix $B(h)$ depending on the
spectral parameter $h$, has the form
$$
B(h)=\left(\begin{array}{ccccc}
0&a_{1}&\cdots&\cdots&-a_{N}h^{-1}\\
-a_{1}&0&\vdots &&\vdots \\
\vdots&\ddots&\ddots&\ddots&\vdots\\
\vdots&&\ddots&\ddots&a_{N-1}\\
a_{N}h&\cdots&\cdots&-a_{N-1}&0
\end{array}\right),
$$
and the rest follows from the general theory. Note that if
$a_j(0)\neq0$, then $a_j(t)\neq0$ for all $t$. Since $A^\top
(h)=A(h^{-1})$, then
$$P(h,z)=\det (A(h)-zI)=P(h^{-1},z).$$
Therefore, the application
\begin{equation}\label{eqn:euler}
\sigma : \mathcal{C}\longrightarrow \mathcal{C},\quad
(h,z)\longmapsto(h^{-1},z),
\end{equation}
is an involution on the spectral curve $\mathcal {C}$. We choose
\begin{eqnarray}
&&A(h)= \left(\begin{array}{ccc}
0&...&a_N\\
\vdots&\ddots&\vdots\\
0&...&0
\end{array}\right)h^{-1}+
\left(\begin{array}{ccccc}
b_1&a_1&&&\\
a_1&b_2&&&\\
&&\ddots&&\\
&&&b_{N-1}&a_{N-1}\\
&&&b_N&a_N
\end{array}\right)\nonumber\\
&&\qquad\qquad+ \left(\begin{array}{ccc}
0&...&0\\
\vdots&\ddots&\vdots\\
a_N&...&0
\end{array}\right)h.\nonumber
\end{eqnarray}
Note that here the matrix $A$ is meromorphic (whereas previously
we considered it to be a polynomial in $h$) but we will see that
we can adopt the theory explained in this section, to this
situation too . We have
$$P(h,z)=-\prod_{j=1}^{N-1}a_j.(h+h^{-1})+z^N+c_1z^{N-1}+\cdots+c_N.$$
Let us assume that $\prod_{j=1}^{N-1}a_j\neq 0$ and pose
\begin{eqnarray}
Q(h,z)\equiv \frac{P(h,z)}{\prod_{j=1}^{N-1}a_j}
&=&h+h^{-1}+\frac{z^N+c_1z^{N-1}+\cdots+c_N}{\prod_{j=1}^{N-1}a_j},\nonumber\\
&=&h+h^{-1}+d_0z^N+d_1z^{N-1}+\cdots+d_N.\nonumber
\end{eqnarray}
In $\mathbb{P}^2(\mathbb{C})$, the affine algebraic curve of
equation $Q(h,z)=0$, is singular at infinity for $n\geq 4$. We
will compute the genus of normalization $\mathcal{C}$ of this
curve. Note that $\mathcal{C}$ is a double covering of
$\mathbb{P}^1(\mathbb{C})$ branched into $2N$ points coinciding
with the fixed points of involution $\sigma$ (3.12), that is,
points where $h =\pm1$. According to the Riemann-Hurwitz formula,
the genus $g$ of the curve $\mathcal{C}$ is
$$g=2\left(g(\mathbb{P}^1(\mathbb{C}))-1\right)+1+\frac{2N}{2}=N-1.$$
Consider the covering $\mathcal{C}\longrightarrow
\mathbb{P}^1(\mathbb{C})$ below and set
$$\frac{1}{z}(\infty)=\mathcal{P}+\mathcal{Q},$$ where $\mathcal{P}$
and $\mathcal{Q}$ are located on two separate sheets. From the
equation $Q(h, z)=0$, the divisor of $h$ is
$$(h)=N\mathcal{P}-N\mathcal{Q}.$$ In that case, the divisor
$\mathcal{D}$ is written
$$\mathcal{D}=N\mathcal{P}+N\mathcal{Q},$$ hence $B\in
H^0(\mathcal{C},\mbox{Hom}(V,V(\mathcal{D}))$. The residue $\rho
(B)\in H^0(\mathcal{C},\mathcal{O}_{\mathcal{D}}(\mathcal{D})$
satisfies the conditions of theorem 7, and consequently the linear
flow is given by the application (2.9). To compute the residue
$\rho(B)$ of $B$, we will determine a set of holomorphic
eigenvectors, using the van Moerbeke-Mumford method described
above. Let us calculate the residue in $\mathcal{Q}$ and the
result will be similarly deduced in $\mathcal{P}$. Let
$\mathcal{E}=\sum_{j=1}^gr_j$ be a general divisor of degree $g$
such that : $$\forall k,
\quad\dim\mathcal{L}(\mathcal{E}+(k-1)\mathcal{P}-k\mathcal{Q})=0.$$
According to the Riemann-Roch theorem,
$$\dim\mathcal{L}(\mathcal{E}+k\mathcal{P}-k\mathcal{Q})\geq 1,$$
hence
$$\forall k,
\quad\dim\mathcal{L}(\mathcal{E}+k\mathcal{P}-k\mathcal{Q})=1.$$
Let $$(f_k)\in
\mathcal{L}(\mathcal{E}+k\mathcal{P}-k\mathcal{Q})=H^0(\mathcal{C},\mathcal{O}_{\mathcal{C}}
(\mathcal{E}+k\mathcal{P}-k\mathcal{Q})), \quad1\leq k\leq N$$ be
a base with $f_N=h$. We can choose a vector $v$ of the following
form $v=(f_1,..., f_N)^\top$, such that $v$ is an eigenvector of
$A$, i.e., $$Av=zv, \quad(h,z)\in \mathbb{C}$$ Hence, $V=h^{-1}v$
is a holomorphic eigenvector. Without restricting generality, we
take $N=3$. The system $Av=zv$, is written explicitly
\begin{eqnarray}
b_1f_1+a_2f_2+a_3&=&zf_1,\nonumber\\
a_1f_1+b_2f_2+a_2h&=&zf_2,\nonumber\\
a_3hf_1+a_2f_2+b_3h&=&zh.\nonumber
\end{eqnarray}
By multiplying each equation of this system by $h^{-1}$,
everything becomes holomorphic except the last equation, i.e.,
$a_3f_1=z+\mbox{Taylor}$. Recall that the residue $\rho(B)$ of $B$
is the section of $\mathcal{O}_{\mathcal{D}}(\mathcal{D})$ induced
by $\lambda$ in the equation (2.7): $Bv=\dot{v}+\lambda v$. In
other words, $$Bv=\varrho(B)v+\mbox{Taylor},$$ and therefore
$$
\left(\begin{array}{c}
\frac{a_1f_2}{h}-\frac{a_3}{h}\\
-\frac{a_1f_1}{h}+a_2\\
a_3f_1-\frac{a_2f_2}{h}
\end{array}\right)=
\left(\begin{array}{c}
0\\
0\\
z
\end{array}\right)+\mbox{Taylor}.$$
We deduce that $\rho(B)=h^{-1}z$, and $\rho(\dot{B})=0$. The same
conclusion holds for the residue in $\mathcal {P}$. Consequently,
the flow in question linearizes on the Jacobian variety of
$\mathcal{C}$.
\end{Exmp}

\section{Algebraically integrable systems}

Consider the system of nonlinear differential equations
\begin{eqnarray}
\frac {dz_{1}}{dt}&=&f_{1}\left(t,
z_{1},...,z_{n}\right),\nonumber\\
&\vdots&\\
\frac{dz_{n}}{dt}&=&f_{n}\left(t, z_{1},...,z_{n}\right),\nonumber
\end{eqnarray}
where $f_{1},...,f_{n}$ are functions of $n+1$ complex variables
$t, z_{1},...,z_{n}$ and which apply a domain of
$\mathbb{C}^{n+1}$ into  $\mathbb{C}$. The Cauchy problem is the
search for a solution $\left(z_{1}(t) ,...,z_{n}(t)\right)$ in a
neighborhood of a point $t_{0}$, satisfying the initial conditions
: $$z_{1}(t_{0})=z_{1}^{0},..., z_{n}(t_{0})=z_{n}^{0}.$$ The
system (4.1) can be written in vector form in $\mathbb{C}^{n}$,
$$\frac {dz}{dt}=f(t,z(t), \quad z=(z_{1},...,z_{n}), \quad
f=(f_{1},...,f_{n}).$$ In this case, the Cauchy problem will be to
determine the solution $z(t)$ such that
$$z(t_{0})=z_{0}=(z_{1}^{0},...,z_{n}^{0}).$$ When the functions
$f_{1}, ...,f_{n}$ are holomorphic in the neighborhood of the
$\left(t_{0},z_{1}^{0},...,z_{n}^{0}\right)$, then the Cauchy
problem admits a holomorphic solution and only one. A question
arises : can the Cauchy problem admits some non-holomorphic
solution in the neighborhood of point
$\left(t_{0},z_{1}^{0},...,z_{n}^{0}\right)?$ When the functions
$f_{1},...,f_ {n}$ are holomorphic, the answer is negative. Other
circumstances may arise for the Cauchy problem concerning the
system of differential equations (4.1), when the holomorphic
hypothesis relative to the functions $f_{1},...,f_{n}$ is no
longer satisfied in the neighborhood of a point. In such a case,
it can be seen that the behavior of the solutions can take on the
most diverse aspects. In general, the singularities of the
solutions are of two types : mobile or fixed, depending on whether
or not they depend on the initial conditions. Important results
have been obtained by Painlev\'{e} [50]. Suppose that the system
(4.1) is written in the form
\begin{eqnarray}
\frac
{dz_{1}}{dt}&=&\frac{P_{1}(t,z_{1},...,z_{n})}{Q_{1}(t,z_{1},...,z_{n})},\nonumber\\
&\vdots&\nonumber\\
\frac{dz_{n}}{dt}&=&\frac{P_{n}(
t,z_{1},...,z_{n})}{Q_{n}(t,z_{1},...,z_{n})},\nonumber
\end{eqnarray}
where
$$P_{k}\left(t,z_{1},...,z_{n}\right) =\sum_{0\leq i_{1},\ldots
,i_{n}\leq p}A_{i_{1},\ldots ,i_{n}}^{\left(k\right) }\left(
t\right) z_{1}^{i_{1}}...z_{n}^{i_{n}},\text{ }1\leq k\leq n,$$
$$ Q_{k}\left(t,z_{1},...,z_{n}\right)=\sum_{0\leq
j_{1},\ldots ,j_{n}\leq q}B_{j_{1},\ldots ,j_{n}}^{\left( k\right)
}(t) z_{1}^{j_{1}}...z_{n}^{j_{n}},\text{ }1\leq k\leq n,$$
polynomials with several indeterminate $z_{1},...,z_{n}$ and
algebraic coefficients in $t$. There are two cases :

$\textbf{(i)}$ the fixed singularities are constituted by four
sets of points. The first is the set of singular points of the
coefficients $A_{i_{1},\ldots ,i_{n}}^{\left(
k\right)}\left(t\right)$, $B_{j_{1},\ldots ,j_{n}}^{\left(
k\right)}\left(t\right)$ intervening in the polynomials
$P_{k}\left(t,z_{1},...,z_{n}\right)$ and
$Q_{k}\left(t,z_{1},...,z_{n}\right)$. In general this set
contains $t=\infty$. The second set consists of the points
$\alpha_{j}$ such that : $Q_{k}\left(t,z_{1},...,z_{n}\right)=0$,
which occurs if all the coefficients
$B_{j_{1},\ldots,j_{n}}^{\left(k\right)}\left(t\right)$ vanish for
$t=\alpha_{j}$. The third is the set of points $\beta_{l}$ such
that for some values $\left(z_{1^{\prime}},...,z_{n^{\prime
}}\right)$ of $\left(z_{1},...,z_{n}\right)$, we have $P_{k}\left(
\beta_{l},z_{1^{\prime }},...,z_{n^{\prime }}\right)=Q_{k}\left(
\beta_{l},z_{1^{\prime }},...,z_{n^{\prime }}\right)=0$. Then the
second members of the above system are presented in the
indeterminate form $\frac{0}{0}$ at the points
$\left(\beta_{l},z_{1^{\prime}},... ,z_{n^{\prime}}\right)$.
Finally, the set of points $\gamma_{n}$ such that there exist
$u_{1},...,u_{n}$, for which
$R_{k}\left(\gamma_{n},u_{1},...,u_{n}\right)=S_{k}\left(
\gamma_{n},u_{1},...,u_{n}\right)=0$, where $R_{k}$ and $S_{k}$
are polynomials in $u_{1},...,u_{n}$ obtained from $P_{k}$ and
$Q_{k}$ by setting $z_{1}=\frac{1}{u_{1}},\ldots,
z_{n}=\frac{1}{u_{n}}$. Each of these sets contains only a finite
number of elements. The system in question has a finite number of
fixed singularities.

$\textbf{(ii)}$ the mobile singularities of solutions of this
system are algebraic mobile singularities: poles and (or)
algebraic critical points. There are no essential singular points
for the solution $\left(z_{1},...,z_{n}\right)$.

Considering the system of differential equations (4.1), we can
find sufficient conditions for the existence and uniqueness of
meromorphic solutions. The existence and uniqueness for the
solution of the Cauchy problem concerning the system (4.1), can be
obtained using the method of indeterminate coefficients. The
solution will be explained in the form of a Laurent series. Such a
solution is formal because we obtain it by performing on various
series, which we assume a priori convergent, various operations
whose validity remains to be justified. The problem of convergence
will therefore arise. The result will therefore be established as
soon as we have verified that these series are convergent. This
will be done using the majorant method [14, 7, 21]. In the
following, we will consider the Cauchy problem concerning the
normal system (4.1) where $f_{1},...,f_{n}$ do not depend
explicitly on $t$, i.e.,
\begin{eqnarray}
\frac {dz_{1}}{dt}&=&f_{1}\left(z_{1},...,z_{n}\right),\nonumber\\
&\vdots&\\
\frac{dz_{n}}{dt}&=&f_{n}\left(z_{1},...,z_{n}\right).\nonumber
\end{eqnarray}
We suppose that $f_{1},...,f_{n}$ are rational functions in
$z_{1},...,z_{n}$ and that the system (4.2) is weight-homogeneous,
i.e., there exist positive integers $l_{1},...,l_{n}$ such that
$$f_{i}(\alpha ^{l_{1}}z_{1},...,\alpha ^{l_{n}}z_{n})=\alpha
^{l_{i}+1}f_{i}(z_{1},...,z_{n}), \quad1\leq i\leq n,$$ for each
non-zero constant $\alpha$. In other words, the system (4.2) is
invariant under the transformation $$t\rightarrow
\alpha^{-1}t,\text{ }z_{1}\rightarrow \alpha^{l_{1}}z_{1},\ldots
,\text{ }z_{n}\rightarrow \alpha^{l_{n}}z_{n}.$$ Note that if the
determinant $\det\left(z_{j}\frac{\partial f_{i}}{\partial
z_{j}}-\delta_{ij}f_{i}\right)_{1\leq i,j\leq n}$, is not
identically zero, then the choice of the numbers $s_{1},...,s_{n}$
is unique. In what follows, we will assume that $t_0=z_0=0$, which
does not affect the generality of the results.

\begin{Theo}
Suppose that
\begin{equation}\label{eqn:euler}
z_{i}=\frac{1}{t^{k_{i}}}\sum_{k=0}^{\infty }z_{i}^{(k)
}t^{k},\quad 1\leq i\leq n, \quad z^{(0)}\neq0
\end{equation}
($k_i\in\mathbb{Z}$, some $k_i>0$) is the formal solution (Laurent
series), obtained by the method of undetermined coefficients of
the weight-homogeneous system (4.2). Then the coefficients
$z_{i}^{(0)}$ satisfy the nonlinear equation
$$
k_{i}z_{i}^{(0)}+f_{i}(z_{1}^{(0)},...,z_{n}^{(0)})=0,
$$
where $1\leq i\leq n$, while $z_{i}^{(1)},z_{i}^{(2)},...$ each
satisfy a system of linear equations of the form
$$
(L-k\mathcal{I})z^{(k)}=\mbox{some polynomial in the }
z^{(j)},\quad 0\leq j\leq k,
$$
where $z^{(k)}=(z_1^{(k)},...,z_n^{(k)})^\top$ and
$$L\equiv\left(\frac{\partial f_{i}}{\partial
z_{j}}(z^{(0)})+\delta_{ij}k_{i}\right)_{1\leq i,j\leq n},$$ is
the Jacobian matrix. Moreover, the formal series (4.3) are
convergent.
\end{Theo}

The series (4.3) is the only meromorphic solution in the sense
that this solution results from the fact that the coefficients
$z_{i}^{(k)}$ are determined unequivocally with the adopted method
of calculation. The result of the previous theorem applies to the
following quasi-homogeneous differential equation of order $n$ :
$$
\frac{d^{n}z}{dt^{n}}=f\left(z,\frac{dz}{dt},...,\frac{d^{n-1}z}{dt^{n-1}}\right).
$$
$f$ being a rational function in
$\displaystyle{z,\frac{dz}{dt},...,\frac{d^{n-1}z}{dt^{n-1}}}$ and
$$z(t_{0})=z_{1}^{0},\frac{dz}{dt}(t_{0})=z_{2}^{0},...,
\frac{d^{n-1}z}{dt^{n-1}}(t_{0})=z_{n}^{0}.$$ Indeed, the
differential equation above reduces to a system of $n$ first order
differential equations by setting
$$
z\left(t\right)=z_{1}\left(t\right) ,\quad \frac{dz}{dt}\left(
t\right)=z_{2}\left(t\right) ,..., \frac{d^{n-1}z}{dt^{n-1}}\left(
t\right)=z_{n}\left(t\right).
$$
We thus obtain
\begin{eqnarray}
\frac{dz_{1}}{dt}&=&z_{2},\nonumber\\
\frac{dz_{2}}{dt}&=&z_{3},\nonumber\\
&\vdots&\nonumber\\
\frac{dz_{n-1}}{dt}&=&z_{n},\nonumber\\
\frac{dz_{n}}{dt}&=&f\left(z_{1},z_{2},...,z_{n}\right).\nonumber
\end{eqnarray}
Such a system constitutes a particular case of the normal system
(4.2).\\

Consider now Hamiltonian dynamical systems of the form
\begin{equation}\label{eqn:euler}
X_{H}:\dot {z}=J\frac{\partial H}{\partial z}\text{ }\equiv
f(z),\text{ }z\in \mathbb{R}^{m},\quad
\left(^.\equiv\frac{d}{dt}\right)
\end{equation}
where $H$\ is the Hamiltonian and $J=J(z)$ is a skew-symmetric
matrix with polynomial entries in $z$, for which the corresponding
Poisson bracket $$\{H_{i},H_{j}\}=\left\langle \frac{\partial
H_i}{\partial z},J\frac{\partial H_j} {\partial z}\right\rangle,$$
satisfies the Jacobi identities.

The system (4.4) with polynomial right hand side will be called
algebraic complete integrable (in abbreviated form : a.c.i.) in
the sense of Adler-van Moerbeke [7, 29, 58] when the following
conditions hold.\\

$\textbf{i)}$ The system admits $n+k$ independent polynomial
invariants $H_{1},...,H_{n+k}$ of which $k$ invariants (Casimir
functions) lead to zero vector fields $$J\frac{\partial
H_{i}}{\partial z}(z)=0, \quad1\leq i\leq k,$$ the $n=(m-k)/2$
remaining ones $H_{k+1}=H$,...,$H_{k+n}$ are in involution (i.e.,
$\left\{H_{i},H_{j}\right\}=0$),which give rise to $n$ commuting
vector fields. For generic $c_{i}$, the invariant manifolds
$\overset{n+k}{\underset{i=1}{\bigcap }}\left\{z\in
\mathbb{R}^{m}:H_{i}=c_{i}\right\}$ are assumed compact and
connected. According to the Arnold-Liouville theorem [6], there
exists a diffeomorphism $$\overset{n+k}{\underset{i=1}{\bigcap
}}\left\{z\in \mathbb{R}^{m}:H_{i}=c_{i}\right\} \longrightarrow
\mathbb{R}^{n}/Lattice,$$ and the solutions of the system (4.4)
are straight lines motions on these real tori.

$\textbf{ii)}$ The invariant manifolds thought of as lying in
$\mathbb{C}^{m}$,
$$\mathcal{A}=\overset{n+k}{\underset{i=1}{\bigcap }}\{z\in
\mathbb{C}^{m}: H_{i}=c_{i}\},$$ are related, for generic $c_i$,
to Abelian varieties $T^n=\mathbb{C}^{n}/Lattice$ (complex
algebraic tori) as follows :
$$\mathcal{A}=T^n\backslash\mathcal{D},$$ where $\mathcal{D}$ is a
divisor (codimension one subvarieties) in $T^n$. In the natural
coordinates $(t_{1},...,t_{n})$ of $T^n$ coming from
$\mathbb{C}^{n}$, the coordinates $z_{i}=z_{i}(t_{1},...,t_{n})$
are meromorphic and $\mathcal{D}$ is the minimal divisor on $T^n$
where the variables $z_i$ blow up. Moreover, the flows (4.4) (run
with complex time) are straight-line motions on $T^n$.

Mumford gave one in his Tata lectures [48], which includes as well
the noncompact. Algebraic means that the torus\ can be defined as
an intersection $$\bigcap_i\{Z\in\mathbb{P}^N : P_{i}(Z)=0\},$$
involving a large number of homogeneous polynomials $P_{i}$.
Condition $i)$ means, in particular, there is an algebraic map
$$(z_{1}(t),...,z_{m}(t))\longmapsto(s_{1}(t),...,s_{n}(t)),$$
making the following sums linear in $t$ :
$$\sum_{i=1}^{n}\int_{s_{i}(0)}^{s_{i}(t)}\omega _{j}=d
_{j}t\text{ },\text{ }1\leq j\leq n,\text{ }d _{j}\in
\mathbb{C},$$ where $\omega _{1},...,\omega _{n}$\ denote
holomorphic differentials on some algebraic curves. If the
Hamiltonian flow (4.4) is a.c.i., it means that the variables
$z_{i}$\ are meromorphic on the torus $T^n$ and by compactness
they must blow up along a codimension one subvariety (a divisor)
$\mathcal{D}\subset T^n$. By the a.c.i. definition, the flow (4.4)
is a straight line motion in $T^n$ and thus it must hit the
divisor $\mathcal{D}$\ in at least one place. Moreover through
every point of $\mathcal{D}$, there is a straight line motion and
therefore a Laurent expansion around that point of intersection.
Hence the differential equations must admit Laurent expansions
which depend on the $n-1$ parameters defining $\mathcal{D}$ and
the $n+k$ constants $c_{i}$ defining the torus $T^n$, the total
count is therefore $$m-1=\dim(\mbox{phase space})-1,$$ parameters.
The fact that algebraic complete integrable systems possess
$(m-1)$-dimensional families of Laurent solutions, was implicitly
used by Kowalewski [32] in her classification of integrable rigid
body motions. Such a necessary condition for algebraic complete
integrability can be formulated as follows [5] :

\begin{Theo}
If the Hamiltonian system (4.4) (with invariant tori not
containing elliptic curves) is algebraic complete integrable, then
each $z_i$ blows up after a finite (complex) time, and for every
$z_i$, there is a family of solutions (4.4) depending on
$\dim(\mbox{phase space})-1=m-1$ free parameters. Moreover, the
system (4.4) possesses families of Laurent solutions depending on
$m-2$, $m-3$,...,$m-n$ free parameters. The coefficients of each
one of these Laurent solutions are rational functions on affine
algebraic varieties of dimensions $m-1$, $m-2$, $m-3$,...,$m-n$.
\end{Theo}

The question is whether this criterion is sufficient. The main
problem will be to complete the affine variety $\mathcal{A}$, into
an Abelian variety. A naive guess would be to take the natural
compactification $\overline{\mathcal{A}}$ of $\mathcal{A}$ by
projectivizing the equations. Indeed, this can never work for a
general reason: an Abelian variety $\widetilde{\mathcal{A}}$ of
dimension bigger or equal than two is never a complete
intersection, that is it can never be described in some projective
space $\Bbb{P}^{n}$ by $n$-dim $\widetilde{\mathcal{A}}$ global
polynomial homogeneous equations. In other words, if $\mathcal{A}$
is to be the affine part of an Abelian variety,
$\overline{\mathcal{A}}$ must have a singularity somewhere along
the locus at infinity. The trajectories of the vector fields (4.4)
hit every point of the singular locus at infinity and ignore the
smooth locus at infinity. In fact, the existence of meromorphic
solutions to the differential equations (4.4) depending on some
free parameters can be used to manufacture the tori, without ever
going through the delicate procedure of blowing up and down.
Information about the tori can then be gathered from the divisor.
More precisely, around the points of hitting, the system of
differential equations (4.4) admit a Laurent expansion solution
depending on $m-1$ free parameters and in order to regularize the
flow at infinity, we use these parameters to blowing up the
variety $\overline{\mathcal{A}}$ along the singular locus at
infinity. The new complex variety obtained in this fashion is
compact, smooth and has commuting vector fields on it; it is
therefore an Abelian variety. The system (4.4) with $k+n$
polynomial invariants has a coherent tree of Laurent solutions,
when it has families of Laurent solutions in $t$, depending on
$n-1$, $n-2$,..., $m-n$ free parameters. Adler and van Moerbeke
[5] have shown that if the system possesses several families of
$(n-1)$-dimensional Laurent solutions (principal Painlev\'{e}
solutions) they must fit together in a coherent way and as we
mentioned above, the system must possess $(n-2)$-,
$(n-3)$-,...dimensional Laurent solutions (lower Painlev\'{e}
solutions), which are the gluing agents of the $(n-1)$-dimensional
family. The gluing occurs via a rational change of coordinates in
which the lower parameter solutions are seen to be genuine limits
of the higher parameter solutions an which in turn appears due to
a remarkable propriety of algebraic complete integrable systems;
they can be put into quadratic form both in the original variables
and in their ratios. As a whole, the full set of Painlev\'{e}
solutions glue together to form a fiber bundle with singular base.
A partial converse to theorem 12, can be formulated as follows [5]
:
\begin{Theo}
If the Hamiltonian system (4.4) satisfies the condition $i)$ in
the definition of algebraic complete integrability and if it
possesses a coherent tree of Laurent solutions, then the system is
algebraic complete integrable and there are no other
$m-1$-dimensional Laurent solutions but those provided by the
coherent set.
\end{Theo}

We assume that the divisor is very ample and in addition
projectively normal (see [17, 42] for definitions when needed).
Consider a point $p\in\mathcal{D}$, a chart $U_j$ around $p$ on
the torus and a function $y_j$ in $\mathcal{L}(\mathcal{D})$
having a pole of maximal order at $p$. Then the vector $(1/y_j,
y_1/y_j,\ldots,y_N/y_j)$ provides a good system of coordinates in
$U_j$. Then taking the derivative with regard to one of the flows
$$\left(\frac{y_i}{y_j}\right)\dot{}=\frac{\dot{y_i}y_j-y_i\dot{y_j}}{y_j^2},\quad
1\leq j\leq N,$$ are finite on $U_j$ as well. Therefore, since
$y_j^2$ has a double pole along $\mathcal{D}$, the numerator must
also have a double pole (at worst), i.e.,
$\dot{y_i}y_j-y_i\dot{y_j}\in \mathcal{L}(2\mathcal{D})$. Hence,
when $\mathcal{D}$ is projectively normal, we have that
$$\left(\frac{y_i}{y_j}\right)\dot{}=\sum_{k,l}a_{k,l}\left(\frac{y_k}{y_j}\right)\left(\frac{y_l}{y_j}\right),$$
i.e., the ratios $y_i/y_j$ form a closed system of coordinates
under differentiation. At the bad points, the concept of
projective normality play an important role: this enables one to
show that $y_i/y_j$ is a bona fide Taylor series starting from
every point in a neighborhood of the point in question. Moreover,
the Laurent solutions provide an effective tool for find the
constants of the motion. For that, just search polynomials $H_i$
of $z$, having the property that evaluated along all the Laurent
solutions $z(t)$ they have no polar part. Indeed, since an
invariant function of the flow does not blow up along a Laurent
solution, the series obtained by substituting the formal solutions
(4.3) into the invariants should, in particular, have no polar
part. The polynomial functions $H_i(z(t))$ being holomorphic and
bounded in every direction of a compact space, (i.e., bounded
along all principle solutions), are thus constant by a Liouville
type of argument. It thus an important ingredient in this argument
to use all the generic solutions. To make these informal arguments
rigorous is an outstanding question of the subject. Assume
Hamiltonian flows to be weight-homogeneous with a weight
$l_{i}\in\mathbb{N}$, going with each variable $z_{i}$. Observe
that then the constants of the motion $H$ can be chosen to be
weight-homogeneous : $$H\left(\alpha
^{l_{1}}z_{1},...,\alpha^{l_{m}}z_{m}\right)=
\alpha^{k}H\left(z_{1},...,z_{m}\right), \quad k\in \mathbb{Z}$$
The study of the algebraic complete integrability of Hamiltonian
systems, includes several passages to prove rigorously. Here we
mention the main passages. We saw that if the flow is
algebraically completely integrable, the differential equations
(4.4) must admits Laurent series solutions (4.3) depending on
$m-1$ free parameters. We must have $k_{i}=l_{i}$ and coefficients
in the series must satisfy at the 0$^{th}$step non-linear
equations,
\begin{equation}\label{eqn:euler}
f_{i}\left(z_{1}^{\left(0\right)},...,z_{m}^{\left( 0\right)
}\right)+g_{i}z_{i}^{\left(0\right) }=0,\text{ }1\leq i\leq m,
\end{equation}
and at the $k$$^{th}$step, linear systems of equations :
\begin{equation}\label{eqn:euler}
\left( L-kI\right) z^{\left( k\right) }=
\left\{\begin{array}{rl} 0&\mbox{ for } k=1\\
\mbox{some polynomial in}& z^{\left( 1\right) },...,z^{\left(
k-1\right)} \mbox{ for } k>1,
\end{array}\right.
\end{equation}
where
$$L=\text{ Jacobian map of }(4.5)=
\text{ }\frac{\partial f}{\partial z}+gI\mid _{z=z^{\left(
0\right)}}.$$ If $m-1$ free parameters are to appear in the
Laurent series, they must either come from the non-linear
equations (4.5) or from the eigenvalue problem (4.6), i.e., $L$
must have at least $m-1$ integer eigenvalues. These are much less
conditions than expected, because of the fact that the homogeneity
$k$ of the constant $H$ must be an eigenvalue of $L$. The formal
series solutions are convergent as a consequence of the majorant
method. Thus, the first step is to show the existence of the
Laurent solutions, which requires an argument precisely every time
$k$ is an integer eigenvalue of $L$ and therefore $L-kI$ is not
invertible. One shows the existence of the remaining constants of
the motion in involution so as to reach the number $n+k$. Then you
have to prove that for given $c_{1},...,c_{m},$ the set
$$\mathcal{D}\equiv \left\{\begin{array}{rl}
&z_{i}(t)=t^{-\nu_{i}}\left(z_{i}^{(0)}+z_{i}^{(1)}t+z_{i}^{(2)}t^{2}+\cdots\right), 1\leq i\leq m\\
&\mbox{Laurent solutions such that}:
H_{j}\left(z_{i}(t)\right)=c_{j}+\mbox{Taylor part}
\end{array}\right\}
$$
defines one or several $n-1$ dimensional algebraic varieties
("Painlev\'{e}" divisor) having the property that
$\overset{n+k}{\underset{i=1}{\bigcap}}\left\{z\in \mathbb{C}^{m}:
H_{i}=c_{i}\right\} \cup \mathcal{D}$, is a smooth compact,
connected variety with $n$ commuting vector fields independent at
every point, i.e., a complex algebraic torus
$\mathbb{C}^{n}/Lattice$. Therefore, the flows $\,J\frac{\partial
H_{k+i}}{\partial z}$ $,...,J\frac{\partial H_{k+n}}{\partial z}$
are straight line motions on this torus (for concrete
applications, see for example [5, 6, 7, 19, 36, 38, 58]). Let's
point out that having computed the space of functions
$\mathcal{L}(\mathcal{D})$ with simple poles at worst along the
expansions, it is often important to compute the space of
functions $\mathcal{L}(k\mathcal{D})$ of functions having $k$-fold
poles at worst along the expansions. These functions play a
crucial role in the study of the procedure for embedding the
invariant tori into projective space. From the divisor
$\mathcal{D}$, a lot of information can be obtained with regard to
the periods and the action-angle variables.

The idea of the Adler-van Moerbeke's proof [4, 40] (or what can be
called the Liouville-Arnold-Adler-van Moerbeke theorem) is closely
related to the geometric spirit of the (real) Arnold-Liouville
theorem [10]. Namely, a compact complex $n$-dimensional variety on
which there exist $n$ holomorphic commuting vector fields which
are independent at every point is analytically isomorphic to a
$n$-dimensional complex torus $\mathbb{C}^{n}/Lattice$ and the
complex flows generated by the vector fields are straight lines on
this complex torus.

\begin{Theo}
Let $\overline{\mathcal{A}}$ be an irreducible variety defined by
an intersection
$$\overline{\mathcal{A}}=\bigcap_i\{Z=(Z_0,Z_1,...,Z_n)\in\mathbb{P}^N(\mathbb{C}):P_i(Z)=0\},$$
involving a large number of homogeneous polynomials $P_i$ with
smooth and irreducible affine part
$$\mathcal{A}=\overline{\mathcal{A}}\cap\{Z_0\neq 0\}.$$ Put
$\overline{\mathcal{A}}\equiv \mathcal{A}\cup \mathcal{D}$, i.e.,
$\mathcal{D}=\overline{\mathcal{A}}\cap \{Z_0=0\}$ and consider
the map $$f:\overline{\mathcal{A}}\longrightarrow
\mathbb{P}^N(\mathbb{C}), \quad Z\longmapsto f(Z).$$ Let
$\widetilde{\mathcal{A}}=f(\overline{\mathcal{A}})=\overline{f(\mathcal{A})}$,
$\mathcal{D}=\mathcal{D}_1\cup...\cup \mathcal{D}_r$, where
$\mathcal{D}_i$ are codimension $1$ subvarieties and
$$\mathcal{S}\equiv f(\mathcal{D})=f(\mathcal{D}_1)\cup...\cup
f(\mathcal{D}_r)\equiv \mathcal{S}_1\cup...\cup \mathcal{S}_r.$$
Assume that :

$(i)$ $f$ maps $\mathcal{A}$ smoothly and 1-1 onto
$f(\mathcal{A})$.

$(ii)$ There exist $n$ holomorphic vector fields $X_1,...,X_n$ on
$\mathcal{A}$ which commute and are independent at every point.
One vector field, say $X_k$ (where $1\leq k\leq n$), extends
holomorphically to a neighborhood of $\mathcal{S}_k$ in the
projective space $\mathbb{P}^N(\mathbb{C})$.

$(iii)$ For all $p\in \mathcal{S}_k$, the integral curve $f(t)\in
\mathbb{P}^N(\mathbb{C})$ of the vector field $X_k$ through
$f(0)=p\in \mathcal{S}_k$ has the property that $$\{f(t) : 0<\mid
t\mid<\varepsilon, t\in \mathbb{C}\}\subset f(\mathcal{A}).$$ This
condition means that the orbits of $X_k$ through $\mathcal{S}_k$
go immediately into the affine part and in particular, the vector
field $X_k$ does not vanish on any point of $\mathcal{S}_k$. Then

$a)$ $\widetilde{A}$ is compact, connected, admits an embedding
into $\mathbb{P}^N(\mathbb{C})$ and is diffeomorphic to a
$n$-dimensional complex torus. The vector fields $X_{1},...,X_{n}$
extend holomorphically and remain independent on
$\widetilde{\mathcal{A}}$.

$b)$ $\widetilde{\mathcal{A}}$ is a K\"{a}hler variety, a Hodge
variety and in particular, $\mathcal{A}$ is the affine part of an
Abelian variety $\widetilde{\mathcal{A}}$.
\end{Theo}

\begin{Exmp}
The periodic $5$-particle Kac-van Moerbeke lattice [23] is given
by the following quadratic vector field
$$\overset{.}{x}_j=x_j(x_{j-1}-x_{j+1}),\quad j=1,...,5$$
o\`{u} $(x_1,...,x_5)\in\mathbb{C}^5$ et $x_j=x_{j+5}$. This
system forms a Hamiltonian vector field for the Poisson structure
$$\{x_j,x_k\}=x_jx_k(\delta_{j,k+1}-\delta_{j+1,k}), \quad 1\leq
j,k\leq 5$$ and admits three independent first integrals
\begin{eqnarray}
H_{1}&=& x_1x_3+x_2x_4+x_3x_5+x_4x_1+x_5x_2,\nonumber\\
H_{2}&=&x_1+x_2+x_3+x_4+x_5,\nonumber\\
H_{3}&=&x_1x_2x_3x_4x_5.\nonumber
\end{eqnarray}
Let's show that this system is algebraically completely
integrable. We easily check that $H_1$ and $H_2$ are involution
while $H_3$ is a Casimir. The system in question is therefore
integrable in the Liouville sense. In addition, it is shown [6]
that the affine variety
$$\bigcap_{j=1}^3\{x\in\mathbb{C}^5 : H_j(x)=c_j\},
(c_1,c_2,c_3)\in\mathbb{C}^3, \quad c_3\neq 0$$ defined by the
intersection of the constants of motion is isomorphic to
$\mbox{Jac}(\mathcal{C})\backslash\mathcal{D}$ where
$$\mathcal{C}=\{(z,w)\in\mathbb{C}^2 : w^2=(z^3-c_1z^2+c_2z)^2-4z\},$$
is a smooth curve of genus $2$ and $\mathcal{D}$ consists of five
copies of $\mathcal{C}$ in the Jacobian variety
$\mbox{Jac}(\mathcal{C})$. The flows generated by $H_1$ et $H_2$
are linearized on $\mbox{Jac}(\mathcal{C})$ and the system in
question is algebraically completely integrable. The reader
interested in the study of this system via various methods can
find further information with more detail in [6] as well as in
[56].
\end{Exmp}

\begin{Exmp} The problem we are going to study now is the generalized periodic
Toda systems. Let $e_0, ...,e_l$ be linearly dependent vectors in
the Euclidean vector space $(\mathbb{R}^{l + 1}, \langle. |.
\rangle)$, $l\geq1$ , such that they are $l$ to $l$ linearly
independent (i.e, for all $j$, the vectors
$e_0,...,\widehat{e_j},..., e_l$ are linearly independent).
Suppose that the non-zero reals $\xi_0,...,\xi_l$ satisfying
$$\sum_{j=0}^l\xi_je_j=0,$$ are non-zero sum; that is,
$$\sum_{j=0}^l\xi_j\neq0.$$ Let $A=(a_{ij})_{0\leq i, j\leq l}$ be
the matrix whose elements are defined by
$$a_{ij}=2\frac{\langle e_i|e_j\rangle}{\langle
e_j|e_j\rangle},\quad 0\leq i,j\leq l.$$ We consider the vector
field $X_A$ on $\mathbb{C}^{2(l+1)}$,
$$
X_A : \left\{\begin{array}{rl}
\overset{.}{x}=x.y\\
\overset{.}{y}=Ax
\end{array}\right.
$$
where $x,y\in\mathbb{C}^{l+1}$ and $x.y=(x_0y_0,...,x_ly_l)$. It
has been shown [7] that if $X_A$ is an integrable vector field of
an irreducibly algebraically completely integrable system, then
$A$ is the Cartan matrix of a possibly twisted affine Lie algebra.
This system was studied by many authors (see [17] and references
therein). Specific detailed results can be found on the technical
paper [6] (and also in [7]) about link between the Toda lattice,
Dynkin diagrams, singularities and Abelian varieties. The periodic
$l+1$ particle Toda lattices are associated to extended Dynkin
diagrams. They have $l+1$ polynomial invariants, as many as there
are dots in the Dynkin diagram and are integrable Hamiltonian
systems. The complex invariant manifold defined by putting these
invariants equals to generic constants completes into an Abelian
variety by gluing on a specific divisor $\mathcal{D}$. The latter
is entirely described by the extended Dynkin diagram : each point
of the diagram corresponds to a component of the divisor and each
subdiagram determines the intersection of the corresponding
divisors. The global geometry of the complex invariant tori
(Abelian varieties), such as polarization, divisor equivalences,
dimension of certain linear systems, etc., is also entirely given
by the extended Dynkin diagram and the linear equivalence between
them is expressed in Lie-theoretic terms. More precisely, the
divisor $\mathcal{D}$ consists of $l+1$ irreducible components
$\mathcal{D}_j$ each associated with a root $\alpha_j$ of the
Dynkin diagram $\Delta$. The intersection of $k$ components
$\mathcal{D}_{{j}_1},...,\mathcal{D}_{{j}_k}$ satisfies the
following relation : the intersection multiplicity of the
intersection of $k$ components of the divisor equals
$\frac{\mbox{order}(W)}{\det(A)}$ where $W$ and $A$ are the Weyl
group and the Cartan matrix going with the sub-Dynkin diagram
$\alpha_{{j}_1},...,\alpha_{{j}_k}$ associated with the $k$
components. The intersection of all the divisors is empty and the
intersection of all divisors but one is a discrete set of points
whose number is explicitly determined. we have the following
expression for the set-theoretical number of points in terms of
the Dynkin diagram
$$\sharp\left(\bigcap_{\beta\neq\alpha}\mathcal{D}_\beta\right)=\frac{p_\alpha}{p_0}
\left(\frac{\mbox{order(Weyl group of the Dynkin diagram
}\Delta\backslash \alpha_0)}{\mbox{order(Weyl group of the Dynkin
diagram }\Delta\backslash \alpha)}\right),$$ where the integers
$p_\alpha$, are given by the null vector of the Cartan matrix
going with $\Delta$. The singularities of the divisor are
canonically associated to semi-simple Dynkin diagrams. The
singularities of each component occur only at the intersections
with other components and their multiplicities at the intersection
with other divisors are expressed in terms of how a corresponding
root is located in the sub-Dynkin diagram determined by this root
and those of the members of the above divisor intersection. The
following inclusion holds for the singular locus
$\mbox{sing}(\mathcal{D}k)$ of $\mathcal{D}_k$ :
$$\mbox{sing}(\mathcal{D}_k)\subseteq\mathcal{D}_k\cap\sum_{\underset{j\neq k}{0\leq j\leq l}}\mathcal{D}_j,\quad
k=0,...,l$$ The multiplicity of the singularity of a particular
component $\mathcal{D}_k$, at its intersection with $m$ other
divisors, i.e.,
$\mbox{sing}(\mathcal{D}_k)\cap(\mathcal{D}_{j_1},...,\mathcal{D}_{j_m})$,
all $j_1,...,j_m\neq k$, is entirely specified by the way the
corresponding root $\alpha_k$ sits in the sub-Dynkin diagram
$\alpha_k, \alpha_{j_1},...,\alpha_{j_m}$. (See [6], for proof of
these results as well as other information).
\end{Exmp}

There are many examples of dynamical systems which have the weak
Painlev\'{e} property that all movable singularities of the
general solution have only a finite number of branches and some
integrable systems appear as coverings of algebraic completely
integrable systems. The invariant varieties are coverings of
Abelian varieties and these systems are called algebraic
completely integrable in the generalized sense. These systems are
Liouville integrable and by the Arnold-Liouville theorem, the
compact connected manifolds invariant by the real flows are tori;
the real parts of complex affine coverings of Abelian varieties.
Most of these systems of differential equations possess solutions
which are Laurent series of $t^{1/n}$ ($t$ being complex time) and
whose coefficients depend rationally on certain algebraic
parameters. In other words, for these systems just replace in the
definition of the complete algebraic integrability above the
condition $\textbf{ii)}$ by the following : $\textbf{iii)}$ the
invariant manifolds $\mathcal{A}$ are related to an $l$-fold cover
$\widetilde{T}^n$ of $T^n$ ramified along a divisor $\mathcal{D}$
in $T^n$ as follows : $\mathcal{A}=\widetilde{T}^n\backslash
\mathcal{D}$.

\begin{Exmp} Let us consider the Ramani-Dorizzi-Grammaticos (RDG) series of integrable potentials [52, 16] :
$$V(x,y)=\sum_{k=0}^{[m/2]}2^{m-2i}\left(\begin{array}{c}
m-i\\
i
\end{array}\right)x^{2i}y^{m-2i},\quad m=1,2,...$$
It can be straightforwardly proven that a Hamiltonian $H$ :
$$H=\frac{1}{2}(p_x^2+p_y^2)+\alpha_mV_m,,\quad m=1,2,...$$
containing $V$ is Liouville integrable, with an additional first
integral :
$$F=p_x(xp_y-yp_x)+\alpha_mx^2V_{m-1},\quad m=1,2,...$$
The study of cases $m=1$ and $m=2$ is easy. The study of other
cases is not obvious. For the case $m=3$, one obtains the
H\'{e}non-Heiles system [19]:
\begin{eqnarray}
\overset{.}{y}_{1}&=&x_{1},\nonumber\\
\overset{.}{y}_{2}&=&x_{2},\nonumber\\
\overset{.}{x}_{1}&=&-Ay_{1}-2y_{1}y_{2},\\
\overset{.}{x}_{2}&=&-16Ay_{2}-y_{1}^{2}-16y_{2}^{2},\nonumber
\end{eqnarray}
corresponding to a generalized H\'{e}non-Heiles Hamiltonian
$$
H=\frac{1}{2}(x_{1}^{2}+x_{2}^{2})+\frac{A}{2}%
(y_{1}^{2}+16y_{2}^{2})+y_{1}^{2}y_{2}+\frac{16}{3}y_{2}^{3},
$$
where $A$ is a constant parameters and $y_{1},y_{2},x_{1},x_{2}$
are canonical coordinates and momenta, respectively. The system
(4.7) can be written in the form
$$
\dot{u}=J\frac{\partial H}{\partial u},\quad
u=(y_{1},y_{2},x_{1},x_{2})^{\top},
$$
where $$\frac{\partial H}{\partial z}=\left(\frac{\partial
H}{\partial y_{1}},\frac{\partial H}{\partial
y_{2}},\frac{\partial H}{\partial x_{1}},\frac{\partial
H}{\partial x_{2}}\right)^{\intercal}, \qquad J=\left(
\begin{array}{ll}
O & I \\
-I & O
\end{array}
\right).$$ The second integral of motion is
$$H_{2}=3x_{1}^{4}+6Ax_{1}^{2}y_{1}^{2}+12x_{1}^{2}y_{1}^{2}y_{2}-4x_{1}x_{2}y_{1}^{3}-4Ay_{1}^{4}y_{2}
-4y_{1}^{4}y_{2}^{2}+3A^{2}y_{1}^{4}-\frac{2}{3}y_{1}^{6}.$$ When
one examines all possible singularities, one finds that it
possible for the variable $y_1$ to contain square root terms of
the type $\sqrt{t}$, which are strictly not allowed by the
Painlev\'{e} test. However, these terms are trivially removed by
introducing some new variables $z_1,\ldots,z_5$, which restores
the Painlev\'{e} property to the system (to see further). And
reasoning as above, we obtain a new algebraically completely
integrable system. The functions $H_1\equiv H$ and $H_{2}$ commute
: $$\left\{H_{1},H_{2}\right\}=\sum_{k=1}^{2}\left(\frac{\partial
H_{1}}{\partial x_{k}}\frac{\partial H_{2}}{\partial
y_{k}}-\frac{\partial H_{1}}{\partial y_{k}}\frac{\partial
H_{2}}{\partial x_{k}}\right)=0.$$ The second flow commuting with
the first is regulated by the equations :
$$\dot{u}=J\frac{\partial H_{2}}{\partial u}.$$ The system (4.7)
admits Laurent solutions in $\sqrt{t}$, depending on three free
parameters : $\alpha$, $\beta$, $\gamma$ and they are explicitly
given as follows
\begin{eqnarray}
y_{1}&=&\frac{\alpha }{\sqrt{t}}+\beta t\sqrt{t}-\frac{\alpha
}{18}t^{2}\sqrt{t}+\frac{\alpha
A_1^{2}}{10}t^{3}\sqrt{t}-\frac{\alpha ^{2}\beta}{18}t^{4}\sqrt{t}+\cdots,\nonumber\\
y_{2}&=&-\frac{3}{8t^{2}}-\frac{A_1}{2}+\frac{\alpha^{2}}{12}t-\frac{2A_1^{2}}{5}t^{2}
+\frac{\alpha\beta}{3}t^{3}-\gamma t^{4}+\cdots,\\
x_{1}&=&-\frac{1}{2}\frac{\alpha}{t\sqrt{t}}+\frac{3}{2}\beta\sqrt{t}-\frac{5}{36}\alpha
t\sqrt{t}+\frac{7}{20}\alpha A_1^{2}t^{2}\sqrt{t}-\frac{1}{4}\alpha^{2}\beta t^{3}\sqrt{t}+\cdots,\nonumber\\
x_{2}&=&\frac{3}{4t^{3}}+\frac{1}{12}\alpha^{2}-\frac{4}{5}A_1^{2}t+\alpha\beta
t^{2}-4\gamma t^{3}+\cdots\nonumber
\end{eqnarray}
These formal series solutions are convergent as a consequence of
the majorant method. By substituting these series in the constants
of the motion $H_{1}=b_{1}$ and $H_{2}=b_{2}$, i.e.,
\begin{eqnarray}
H_{1}&=&\frac{1}{9}\alpha ^{2}-\frac{21}{4}\gamma
+\frac{13}{288}\alpha ^{4}+\frac{4}{3}A^{3}=b_1,\nonumber\\
H_{2}&=&-144\alpha \beta ^{3}+\frac{294}{5}\alpha ^{3}\beta
A^{2}+\frac{8}{9}\alpha ^{6}-33\gamma \alpha ^{4}=b_2,\nonumber
\end{eqnarray}
one eliminates the parameter $\gamma$ linearly, leading to an
equation connecting the two remaining parameters $\alpha$ and
$\beta$ :
$$
144\alpha \beta ^{3}-\frac{294A^{2}}{5}\alpha ^{3}\beta
+\allowbreak \frac{143}{504}\alpha ^{8}-\frac{4}{21}\alpha
^{6}+\frac{44}{21}\left( 4A^{3}-3b_{1}\right) \alpha ^{4}+b_{2}=0.
$$
which is nothing but the equation of an algebraic curve
$\mathcal{D}$ along which the $u(t)
\equiv(y_{1}(t),y_{2}(t),x_{1}(t),x_{2}(t))$ blow up. To be more
precise $\mathcal{D}$ is the closure of the continuous components
of
$$
\left\{\text{Laurent series solutions }u(t)\text{ such that }
H_{k}(u(t))=b_{k},\text{ }1\leq k\leq 2\right\},
$$
i.e.,
$$
\mathcal{D}==t^{0}-\text{coefficient of }\left\{u\in
\Bbb{C}^{4}:\text{ }H_{1}(u(t))=b_{1}\right\}\cap\left\{u\in
\Bbb{C}^{4}:H_{2}(u(t))=b_{2}\right\} .
$$
The invariant variety
\begin{equation}\label{eqn:euler}
\mathcal{A}=\bigcap_{k=1}^{2}\{z\in\mathbb{C}^4 : H_k(z)=b_k\},
\end{equation}
is a smooth affine surface for generic
$(b_{1},b_{2})\in\Bbb{C}^{2}$. The Laurent solutions restricted to
the surface $\mathcal{A}$ are parameterized by the curve
$\mathcal{D}$. We show that the system (4.7) is part of a new
system of differential equations in five unknowns having two cubic
and one quartic invariants (constants of motion). By inspection of
the expansions (4.8), we look for polynomials in
$(y_1,y_2,x_1,x_2)$ without fractional exponents. Let
\begin{equation}\label{eqn:euler}
\varphi : \mathcal{A}\longrightarrow \mathbb{C}^5,\quad
(y_1,y_2,x_1,x_2)\longmapsto (z_1,z_2,z_3,z_4,z_5), \end{equation}
be a morphism on the affine variety $\mathcal{A}$(4.9) where
$z_1,\ldots,z_5$ are defined as
$$z_1=y_{1}^{2},\quad z_2=y_{2},\quad z_3=x_{2},\quad
z_4=y_{1}x_{1},\quad z_5=3x_{1}^{2}+2y_{1}^{2}y_{2}.$$ Using the
two first integrals $H_1$, $H_2$ and differential equations (4.7),
we obtain a system of differential equations in five unknowns,
\begin{eqnarray}
\dot{z}_1&=&2z_4,\nonumber\\
\dot{z}_2&=&z_3,\nonumber\\
\dot{z}_3&=&-z_{1}-16A_1z_{2}-16z_{2}^{2},\\
\dot{z}_4&=&-A_1z_{1}+\frac{1}{3}z_{5}-\frac{8}{3}z_{1}z_{2},\nonumber\\
\dot{z}_5&=&-6A_1z_{4}+2z_{1}z_{3}-8z_{2}z_{4},\nonumber
\end{eqnarray}
having two cubic and one quartic invariants (constants of motion),
\begin{eqnarray}
F_1&=&\frac{1}{2}A_1z_{1}+\frac{1}{6}z_{5}+8A_1z_{2}^{2}+\frac{1}{2}z_{3}^{2}
+\frac{2}{3}z_{1}z_{2}+\frac{16}{3}z_{2}^{3},\nonumber\\
F_2&=&9A_1^{2}z_{1}^{2}+z_{5}^{2}+6A_1z_{1}z_{5}-2z_{1}^{3}-24A_1z_{1}^{2}z_{2}
-\allowbreak 12z_{1}z_{3}z_{4}+24z_{2}z_{4}^{2}-16z_{1}^{2}z_{2}^{2},\nonumber\\
F_3&=&z_{1}z_{5}-3z_{4}^{2}-2z_{1}^{2}z_{2}.\nonumber
\end{eqnarray}
This new system is completely integrable and can be written as
$$\dot{z}=J\frac{\partial H}{\partial z},
\quad z=(z_{1},z_{2},z_{3},z_{4},z_{5})^\top,$$ where $H=F_{1}$.
The Hamiltonian structure is defined by the Poisson bracket
$$\{F,H\}=\left\langle \frac{\partial F}{\partial z},
J\frac{\partial H}{\partial
z}\right\rangle=\sum_{k,l=1}^{5}J_{kl}\frac{\partial F}{\partial
z_{k}}\frac{\partial H}{\partial z_{l}},$$ where $$\frac{\partial
H}{\partial z}=\left(\frac{\partial H}{\partial
z_{1}},\frac{\partial H}{\partial z_{2}},\frac{\partial
H}{\partial z_{3}},\frac{\partial H}{\partial
z_{4}},\frac{\partial H}{\partial z_{5}}\right)^\top,$$ and
$$J=\left(\begin{array}{ccccc}
0&0&0&2z_1&12z_4\\
0&0&1&0&0\\
0&-1&0&0&-2z_1\\
-2z_1&0&0&0&-8z_{1}z_{2}+2z_{5}\\
-12z_4&0&2z_1&8z_{1}z_{2}-2z_{5}&0
\end{array}\right),$$
is a skew-symmetric matrix for which the corresponding Poisson
bracket satisfies the Jacobi identities.  The second flow
commuting with the first is regulated by the equations
$$\dot{z}=J\frac{\partial F_{2}}{\partial z},
\quad z=(z_{1},z_{2},z_{3},z_{4},z_{5})^\top.$$ These vector
fields are in involution, i.e., $\{F_1,F_2\}=0$, and the remaining
one is Casimir, i.e., $J\frac{\partial F_{3}}{\partial z}=0$. The
invariant variety
\begin{equation}\label{eqn:euler}
\mathcal{B}=\bigcap_{k=1}^{3}\{z\in\mathbb{C}^5:F_k(z)=c_k\},
\end{equation}
is a smooth affine surface for generic values of $c_{1},c_{2}$,
$c_{3}$. The system (4.11) possesses Laurent series solutions
which depend on four free parameters. These meromorphic solutions
restricted to the surface $\mathcal{B}$(4.12) can be read off from
(4.8) and the change of variable (4.10). Following the method
mentioned previously, one find the compactification of
$\mathcal{B}$ into an Abelian surface $\widetilde{B}$, the system
of differential equations (4.11) is algebraic complete integrable
and the corresponding flows evolve on $\widetilde{\mathcal{B}}$.
We show that the invariant surface $\mathcal{A}$(4.9) can be
completed as a cyclic double cover $\overline{\mathcal{A}}$ of an
Abelian surface $\widetilde{\mathcal{B}}$. The system (4.7) is
algebraic complete integrable in the generalized sense. Moreover,
$\overline{\mathcal{A}}$ is smooth except at the point lying over
the singularity of type $A_3$ and the resolution
$\widetilde{\mathcal{A}}$ of $\overline{\mathcal{A}}$ is a surface
of general type. We have shown that the morphism $\varphi$(4.10)
maps the vector field (4.7) into an algebraic completely
integrable system (4.11) in five unknowns and the affine variety
$\mathcal{A}$(4.9) onto the affine part $\mathcal{B}$(4.12) of an
Abelian variety $\widetilde{\mathcal{B}}$. This explains (among
other) why the asymptotic solutions to the differential equations
(4.7) contain fractional powers. All this is summarized as follows
[41]:

\begin{Theo}
The system (4.7) admits Laurent solutions with fractional powers
(i.e., contain square root terms of the type $\sqrt{t}$ which are
strictly not allowed by the Painlev\'{e} test) depending on three
free parameters and is algebraic complete integrable in the
generalized sense. The morphism $\varphi$ (4.10) (which restores
the Painlev\'{e} property) maps this system into a new algebraic
completely integrable system (4.11) in five unknowns.
\end{Theo}
The case $m=4$, corresponds the Ramani Dorizzi Grammaticos (RDG)
system [52, 15],
\begin{eqnarray}
\ddot{q}_1-q_{1}\left(q_{1}^{2}+3q_{2}^{2}\right)&=&0,\\
\ddot{q}_2-q_{2}\left( 3q_{1}^{2}+8q_{2}^{2}\right)&=&0,\nonumber
\end{eqnarray}
corresponding to the Hamiltonian
\begin{equation}\label{eqn:euler}
H_1=\frac{1}{2}(p_{1}^{2}+p_{2}^{2})-\frac{3}{2}q_{1}^{2}q_{2}^{2}
-\frac{1}{4}q_{1}^{4}-2q_{2}^{4}.
\end{equation}
This system is integrable in the sense of Liouville, the second
first integral (of degree $8$) being
\begin{equation}\label{eqn:euler}
H_{2}=p_{1}^{4}-6q_{1}^{2}q_{2}^{2}p_{1}^{2}+q_{1}^{4}q_{2}^{4}
-q_{1}^{4}p_{1}^{2}+q_{1}^{6}q_{2}^{2}+\allowbreak
4q_{1}^{3}q_{2}p_{1}p_{2}-q_{1}^{4}p_{2}^{2}+\frac{1}{4}q_{1}^{8}.
\end{equation}
The first integrals $H_1$ and $H_2$ are in involution, i.e.,
$\left\{H_{1},H_{2}\right\}=0$. The system (4.13) is
weight-homogeneous with $q_1, q_2$ having weight $1$ and $p_1,
p_2$ weight $2$, so that $H_1$(4.14) and $H_2$(4.15) have weight
$4$ and $8$ respectively. When one examines all possible
singularities, one finds that it possible for the variable $q_1$
to contain square root terms of the type $\sqrt{t}$, which are
strictly not allowed by the Painlev\'{e} test. However, we will
see later that these terms are trivially removed by introducing
the variables $z_1,\ldots,z_5$ (4.20) which restores the
Painlev\'{e} property to the system. Let $\mathcal{B}$ be the
affine variety defined by
\begin{equation}\label{eqn:euler}
\mathcal{B}=\bigcap_{k=1}^{2}\{z\in\mathbb{C}^4:H_k(z)=b_k\},\end{equation}
for generic $(b_{1},b_{2}) \in \mathbb{C}^{2}$. The system (4.13)
possesses $3$-dimensional family of Laurent solutions (principal
balances) depending on three free parameters $u, v$ and $w$. There
are precisely two such families, labelled by $\varepsilon=\pm 1$,
and they are explicitly given as follows
\begin{eqnarray}
q_{1}&=&\frac{1}{\sqrt{t}}(u-\frac{1}{4}u^{3}t+vt^{2}
-\frac{5}{128}u^{7}t^{3}+\frac{1}{8}u(\frac{3}{4}u^{3}v
-\frac{7}{256}u^{8}+\allowbreak 3\varepsilon w)t^{4}+\cdots),\nonumber\\
q_{2}&=&\frac{1}{t}(\frac{1}{2}\varepsilon-\frac{1}{4}\varepsilon
u^{2}t +\frac{1}{8}\varepsilon u^{4}t^{2}+\frac{1}{4}\varepsilon
u(\frac{1}{32}u^{5}
-3v) t^{3}+\allowbreak wt^{4}+\cdots),\\
p_{1}&=&\frac{1}{2t\sqrt{t}}(-u-\frac{1}{4}u^{3}t+3vt^{2}
-\frac{25}{128}t^{3}u^{7}+\nonumber\\
&&\qquad \qquad \qquad\frac{7}{8}u(\frac{3}{4}u^{3}v
-\frac{7}{256}u^{8}+3\varepsilon w)t^{4}+\cdots),\nonumber\\
p_{2}&=&\frac{1}{t^{2}}(-\frac{1}{2}\varepsilon
+\frac{1}{8}\varepsilon u^{4}t^{2}+\frac{1}{2}\varepsilon u(
\frac{1}{32}u^{5}-3v)t^{3}+3wt^{4}+\cdots).\nonumber
\end{eqnarray}
These formal series solutions are convergent as a consequence of
the majorant method. By substituting these series in the constants
of the motion $H_{1}=b_{1}$ and $H_{2}=b_{2}$, one eliminates the
parameter $w$ linearly, leading to an equation connecting the two
remaining parameters $u$ and $v$ :
\begin{eqnarray}
\Gamma :
&&\frac{65}{4}uv^{3}+\frac{93}{64}u^{6}v^{2}+\frac{3}{8192}\left(
-9829u^{8}+26112H_{1}\right)
u^{3}v\\&&-\frac{10299}{65536}u^{16}-\allowbreak
\frac{123}{256}H_{1}u^{8}+H_{2}+\frac{15362\,98731}{52}=0.\nonumber
\end{eqnarray}
According to Hurwitz's formula, this defines a Riemann surface
$\Gamma$ of genus $16$. The Laurent solutions restricted to the
affine surface $\mathcal{B}$(4.16) are thus parameterized by two
copies $\Gamma_{-1}$ and $\Gamma_{1}$ of the same Riemann surface
$\Gamma$. Let
\begin{equation}\label{eqn:euler}
\varphi : \mathcal{B}\longrightarrow \mathbb{C}^5,\quad
(q_1,q_2,p_1,p_2)\longmapsto (z_1,z_2,z_3,z_4,z_5),
\end{equation}
be the morphism defined on the affine variety $\mathcal{B}$(4.16)
by
\begin{equation}\label{eqn:euler}
z_1=q_{1}^{2},\quad z_2=q_{2},\quad z_3=p_{2},\quad
z_4=q_{1}p_{1},\quad z_5=p_{1}^{2}-q_{1}^{2}q_{2}^{2}.
\end{equation}
These variables are easily obtained by simple inspection of the
series (4.17). By using the variables (4.20) and differential
equations (4.13), one obtains
\begin{eqnarray}
\dot{z}_1&=&2z_4,\nonumber\\
\dot{z}_2&=&z_3,\nonumber\\
\dot{z}_3&=&z_2(3z_1+8z_2^2),\\
\dot{z}_4&=&z_{1}^{2}+4z_{1}z_{2}^{2}+z_{5},\nonumber\\
\dot{z}_5&=&2z_{1}z_{4}+4z_2^2z_4-2z_{1}z_{2}z_{3}.\nonumber
\end{eqnarray}
This new system on $\mathbb {C}^5$ admits the following three
first integrals
\begin{eqnarray}
F_1&=&\frac{1}{2}z_{5}-z_{1}z_{2}^{2}+\frac{1}{2}z_{3}^{2}
-\frac{1}{4}z_{1}^{2}-2z_{2}^{4},\nonumber\\
F_2&=&z_5^2-z_1^2z_5+4z_1z_2z_3z_4-z_1^2z_3^2+\frac{1}{4}z_1^4-4z_2^2z_4^2 ,\\
F_3&=&z_{1}z_{5}+z_{1}^{2}z_{2}^{2}-z_{4}^{2}.\nonumber
\end{eqnarray}
The first integrals $F_1$ and $F_2$ are in involution , while
$F_3$ is trivial (Casimir function). The invariant variety
$\mathcal{A}$ defined by
\begin{equation}\label{eqn:euler}
\mathcal{A}=\bigcap_{k=1}^{2}\{z:F_k(z)=c_k\}\subset\mathbb{C}^5,
\end{equation}
is a smooth affine surface for generic values of
$(c_{1},c_{2},c_{3}) \in \mathbb{C}^{3}$. The system (4.21) is
completely integrable and possesses Laurent series solutions which
depend on four free parameters $\alpha, \beta, \gamma$ et $\theta$
:
\begin{eqnarray}
z_{1}&=&\frac{1}{t}\alpha -\frac{1}{2}\alpha ^{2}+\beta t
-\frac{1}{16}\alpha \left( \alpha ^{3}+4\beta \right) t^{2}+\gamma t^{3}+\cdots, \nonumber\\
z_{2}&=&\frac{1}{2t}\varepsilon -\frac{1}{4}\varepsilon \alpha
+\frac{1}{8}\varepsilon \alpha ^{2}t
-\frac{1}{32}\varepsilon \left(-\alpha ^{3}+12\beta \right) t^{2}+\theta t^{3}+\cdots, \nonumber\\
z_{3}&=&-\frac{1}{2t^{2}}\varepsilon +\frac{1}{8}\varepsilon
\alpha ^{2}\allowbreak
-\frac{1}{16}\varepsilon \left( -\alpha ^{3}+12\beta \right) t+3\theta t^{2}+\cdots,\\
z_{4}&=&-\frac{1}{2t^{2}}\alpha +\frac{1}{2}\beta
-\frac{1}{16}\alpha \left( \alpha ^{3}+4\beta \right) t+\frac{3}{2}\gamma t^{2}+\cdots, \nonumber\\
z_{5}&=&\frac{1}{2t^{2}}\alpha ^{2}-\frac{1}{4t}\left( \alpha
^{3}+4\beta \right) +\allowbreak \frac{1}{4}\alpha \left( \alpha
^{3}+2\beta \right) -\left( \alpha ^{2}\beta -2\gamma
+4\varepsilon \theta \alpha \right)t +\cdots, \nonumber
\end{eqnarray}
where $\varepsilon=\pm 1$. The convergence of these series is
guaranteed by the majorant method. Substituting these developments
in equations (4.22), one obtains three polynomial relations
between $\alpha, \beta, \gamma$ and $\theta$. Eliminating $\gamma$
and $\theta$ from these equations, leads to an equation connecting
the two remaining parameters $\alpha$ and $\beta$ :
\begin{eqnarray}
\mathcal{C}: &&64\beta ^{3}-16\alpha ^{3}\beta ^{2}-4\left( \alpha
^{6}-32\alpha ^{2}c_{1}-16c_{3}\right) \beta \\
&&+\alpha \left( 32c_{2}-32\alpha ^{4}c_{1}+\alpha ^{8}-16\alpha
^{2}c_{3}\right)=0.\nonumber
\end{eqnarray}
The Laurent solutions restricted to the surface $\mathcal{A}$
(4.23) are thus parameterized by two copies $\mathcal{C}_{-1}$ and
$\mathcal{C}_{1}$ of the same Riemann surface $\mathcal{C}$
(4.25). According to the Riemann-Hurwitz formula, the genus of
$\mathcal{C}$ is $7$. We embed these curves in a hyperplane of
$\mathbb{P}^{15}(\mathbb{C})$ using the sixteen functions :
$$1,\quad z_{1},\quad z_{2},\quad 2z_{5}-z_1^2,\quad z_{3}+2\varepsilon
z_2^2,\quad z_{4}+\varepsilon z_1 z_2,\quad W(f_{1},f_{2}),$$$$
f_{1}(f_{1}+2\varepsilon f_{4}), \quad f_{2}(f_{1}+2\varepsilon
f_{4}), \quad z_{4}(f_{3}+2\varepsilon f_{6}), \quad
z_{5}(f_{3}+2\varepsilon f_{6}),$$$$f_{5}(f_{1}+2\varepsilon
f_{4}), \quad f_{1}f_{2}(f_{3}+2\varepsilon f_{6}),\quad
f_4f_5+W(f_{1},f_{4}),$$$$ W(f_{1},f_{3})+2\varepsilon
W(f_{1},f_{6}), \quad f_3-2z_5+4f_4^2,$$ where
$W(s_j,s_k)\equiv\dot s_j s_k-s_j\dot s_k$ (Wronskian) and we show
that these curves have two points in common in which
$\mathcal{C}_1$ is tangent to $\mathcal{C}_{-1}$. The system
(4.13) is algebraic complete integrable in the generalized sense.
The invariant surface $\mathcal{B}$(4.16) can be completed as a
cyclic double cover $\overline{\mathcal{B}}$ of the Abelian
surface $\widetilde{\mathcal{A}}$, ramified along the divisor
$\mathcal{C}_{1}+\mathcal{C}_{-1}$. Moreover,
$\overline{\mathcal{B}}$ is smooth except at the point lying over
the singularity (of type $A_3$) of
$\mathcal{C}_{1}+\mathcal{C}_{-1}$ (double points of intersection
of the curves $\mathcal{C}_1$ and $\mathcal{C}_{-1}$) and the
resolution $\widetilde{\mathcal{B}}$ of $\overline{\mathcal{B}}$
is a surface of general type. We shall resume the proof of these
results. Observe that the morphism $\varphi$(4.19) is an
unramified cover. The Riemann surface $\Gamma$(4.18) play an
important role in the construction of a compactification
$\overline{\mathcal{B}}$ of $\mathcal{B}$. Let us denote by $G$ a
cyclic group of two elements $\{-1,1\}$ on
$$V_\varepsilon^j=U_\varepsilon^j \times \{\tau \in \mathbb{C} :
0<|\tau|<\delta\},$$ where $\tau=t^{1/2}$ and $ U_\varepsilon^j$
is an affine chart of $\Gamma_\varepsilon$ for which the Laurent
solutions (4.24) are defined. The action of $G$ is defined by
$$(-1)\circ (u,v,\tau)=(-u,-v,-\tau),$$ and is without fixed points
in $V_\varepsilon^j$. So we can identify the quotient
$V_\varepsilon^j / G$ with the image of the smooth map
$h_\varepsilon^j :V_\varepsilon ^j\longrightarrow \mathcal{B}$
defined by the expansions (4.24). We have
$$(-1,1).(u,v,\tau)=(-u,-v,\tau),$$ and
$$(1,-1).(u,v,\tau)=(u,v,-\tau),$$ i.e., $G\times G$ acts separately
on each coordinate. Thus, identifying $V_\varepsilon ^j/G^2$ with
the image of $\varphi\circ h_\varepsilon ^j$ in $\mathcal{A}$.
Note that $\mathcal{B}_\varepsilon ^j=V_\varepsilon ^j/G$ is
smooth (except for a finite number of points) and the coherence of
the $\mathcal{B}_\varepsilon^j$ follows from the coherence of
$V_\varepsilon^j$ and the action of $G$. Now by taking
$\mathcal{B}$ and by gluing on various varieties
$\mathcal{B}_\varepsilon^j\backslash \{\mbox{some points}\}$, we
obtain a smooth complex manifold $\widehat{\mathcal{B}}$ which is
a double cover of the Abelian variety $\widetilde{\mathcal{A}}$
ramified along $\mathcal{C}_1+\mathcal{C}_{-1}$, and therefore can
be completed to an algebraic cyclic cover of
$\widetilde{\mathcal{A}}$. To see what happens to the missing
points, we must investigate the image of $\Gamma \times\{0\}$ in
$\cup \mathcal{B}_\varepsilon^j$. The quotient $\Gamma
\times\{0\}/G$ is birationally equivalent to the Riemann surface
$\Upsilon$ of genus $7$ :
$$\Upsilon : \frac{65}{4}y^{3}+\frac{93}{64}x^{3}y^{2}+\frac{3}{8192}
\left(-9829x^{4}+26112b_{1}\right) x^{2}y\\$$$$+ x\left(
-\frac{10299}{65536}x^{8}-\allowbreak
\frac{123}{256}b_{1}x^{4}+b_{2}+
\frac{15362\,98731}{52}\right)=0,$$ where $y=uv, x=u^2$. The
Riemann surface $\Upsilon$ is birationally equivalent to
$\mathcal{C}$. The only points of $\Upsilon$ fixed under
$(u,v)\longmapsto (-u,-v)$ are the points at $\infty$, which
correspond to the ramification points of the map $$\Gamma
\times\{0\}\overset{2-1}{\longrightarrow }\Upsilon
:(u,v)\longmapsto(x,y),$$ and coincides with the points at
$\infty$ of the Riemann surface $\mathcal{C}$. Then the variety
$\widehat{\mathcal{B}}$ constructed above is birationally
equivalent to the compactification $\overline{\mathcal{B}}$ of the
generic invariant surface $\mathcal{B}$. So
$\overline{\mathcal{B}}$ is a cyclic double cover of the Abelian
surface $\widetilde{\mathcal{A}}$ ramified along the divisor
$\mathcal{C}_1+\mathcal{C}_{-1}$, where $\mathcal{C}_1$ and
$\mathcal{C}_{-1}$ have two points in commune at which they are
tangent to each other. The system (4.13) is algebraic complete
integrable in the generalized sense. Moreover,
$\overline{\mathcal{B}}$ is smooth except at the point lying over
the singularity (of type $A_3$) of
$\mathcal{C}_1+\mathcal{C}_{-1}$. In term of an appropriate local
holomorphic coordinate system $(X,Y,Z),$ the local analytic
equation about this singularity is $X^4+Y^2+Z^2=0$. Let
$\widetilde{\mathcal{B}}$ be the resolution of singularities of
$\overline{\mathcal{B}},$ $\mathcal{X}(\widetilde{\mathcal{B}})$
be the Euler characteristic of $\widetilde{\mathcal{B}}$ and
$p_g(\widetilde{\mathcal{B}})$ the geometric genus of
$\widetilde{\mathcal{B}}$. Then $\widetilde{\mathcal{B}}$ is a
surface of general type with invariants:
$\mathcal{X}(\widetilde{\mathcal{B}})=1$ and
$p_g(\widetilde{\mathcal{B}})=2$. In summary we have [37],

\begin{Theo}
The system (4.13) admits Laurent solutions ,
$$(q_1,q_2,p_1,p_2)=(t^{-1/2},t^{-1},t^{-3/2},t^{-2}) \times
\mbox{ a Taylor series in }t,$$ depending on three free parameters
: $u, v $ and $w$. These solutions restricted to the surface
$\mathcal{B}$(4.16) are parameterized by two copies $\Gamma_{1}$
and $\Gamma_{-1}$ of the Riemann surface $\Gamma$(4.18) of genus
$16$. This system (4.13) is algebraic complete integrable in the
generalized sense and extends to a new system (4.29) of five
differential equations algebraically completely integrable with
three quartics invariants (4.22). Generically, the invariant
manifold $\mathcal{A}$(4.23) defined by the intersection of these
quartics form the affine part of an Abelian surface
$\widetilde{\mathcal{A}}$. The reduced divisor at infinity
$$\widetilde{\mathcal{A}}\setminus
\mathcal{A}=\mathcal{C}_1+\mathcal{C}_{-1},$$ is very ample and
consists of two components $\mathcal{C}_1$ and $\mathcal{C}_{-1}$
of a genus $7$ curve $\mathcal{C}$(4.25). In addition, the
invariant surface $\mathcal{B}$ can be completed as a cyclic
double cover $\overline{\mathcal{B}}$ of the Abelian surface
$\widetilde{\mathcal{A}}$, ramified along the divisor
$\mathcal{C}_{1}+\mathcal{C}_{-1}$. Moreover, $\overline{B}$ is
smooth except at the point lying over the singularity (of type
$A_3$) of $\mathcal{C}_{1}+\mathcal{C}_{-1}$ and the resolution
$\widetilde{B}$ of $\overline{B}$ is a surface of general type
with invariants : $\mathcal{X}(\widetilde{B})=1$ and
$p_g(\widetilde{B})=2$.
\end{Theo}

Consider on the Abelian variety $\widetilde{\mathcal{A}}$ the
holomorphic $1$-forms $dt_1$ and $dt_2$ defined by
$dt_i(X_{F_j})=\delta_{ij}$, where $X_{F_1}$ and $X_{F_2}$ are the
vector fields generated respectively by $F_1$ and $F_2$. Taking
the differentials of $\zeta=1/z_2$ and $\xi=\frac{z_1}{z_2}$
viewed as functions of $t_1$ and $t_2$, using the vector fields
and the Laurent series (4.24) and solving linearly for $dt_1$ and
$dt_2$, we obtain the holomorphic differentials
\begin{eqnarray}
\omega_1&=&dt_{1}|_{\mathcal{C}_{\varepsilon}}=\frac{1}{\Delta}(\frac{\partial
\xi}{\partial t_2}d\zeta-\frac{\partial \zeta}{\partial
t_2}d\xi)|_{\mathcal{C}_{\varepsilon}} =
\frac{8}{\alpha \left(-4\beta +\alpha ^{3}\right) }d\alpha ,\nonumber\\
\omega_2&=&dt_{2}|_{\mathcal{C}_{\varepsilon}}=\frac{1}{\Delta}(\frac{-\partial
\xi}{\partial t_1}d\zeta-\frac{\partial \zeta}{\partial
t_1}d\xi)|_{\mathcal{C}_{\varepsilon}}=\frac{2}{\left(-4\beta
+\alpha^{3}\right) ^{2}}d\alpha ,\nonumber
\end{eqnarray}
with $$\Delta \equiv \frac{\partial \zeta}{\partial
t_1}\frac{\partial \xi}{\partial t_2}-\frac{\partial
\zeta}{\partial t_2}\frac{\partial \xi}{\partial t_1}.$$ The
zeroes of $\omega_2$ provide the points of tangency of the vector
field $X_{F_{1}}$ to $\mathcal{C}_\varepsilon$. We have
$$\frac{\omega_1}{\omega_2}=\frac{4}{\alpha
}\left(-4\beta+\alpha^{3}\right),$$ and $X_{F_{1}}$ is tangent to
$\mathcal{H}_\varepsilon$ at the point covering $\alpha=\infty$.
Note that the reflection $\sigma$ on the affine variety
$\mathcal{A}$ amounts to the flip $$\sigma
:(z_1,z_2,z_3,z_4,z_5)\longmapsto (z_1,-z_2,z_3,-z_4,z_5),$$
changing the direction of the commuting vector fields. It can be
extended to the (-Id)-involution about the origin of
$\mathbb{C}^2$ to the time flip $$(t_1,t_2)\longmapsto
(-t_1,-t_2),$$ on $\widetilde{\mathcal{A}}$, where $t_{1}$ and
$t_{2}$ are the time coordinates of each of the flows
$X_{{F}_{1}}$ and $X_{{F}_{2}}$. The involution $\sigma $ acts on
the parameters of the Laurent solution (3.24) as follows $$\sigma
:(t,\alpha,\beta,\gamma,\theta)\longmapsto
(-t,-\alpha,-\beta,-\gamma,\theta),$$ interchanges the Riemann
surfaces $\mathcal{C}_{\varepsilon}$ and the linear space
$\mathcal{L}$ can be split into a direct sum of even and odd
functions. Geometrically, this involution interchanges
$\mathcal{C}_{1}$ and $\mathcal{C}_{-1}$, i.e.,
$\mathcal{C}_{-1}=\sigma \mathcal{C}_{1}$.

\textbf{Remark :} However, the case $m=5$, corresponds to a system
with an Hamiltonian of the form
$$H=\frac{1}{2}(p_x^2+p_y^2)+y^5+x^2y^3+\frac{3}{16}x^4y.$$ The
corresponding Hamiltonian system admits a second first integral :
$$F=-p_x^2y+p_xp_yx-\frac{1}{2}x^2y^4+\frac{3}{8}x^4y^2+\frac{1}{32}x^6,$$
and admits three $3$-dimensional families solutions $x$, $y$,
which are Laurent series of $t^{1/3}$  : $x=at^{-\frac{1}{3}}$,
$x=bt^{-\frac{2}{3}}$, $b^3=-\frac{2}{9}$, but for which there are
no polynomial $P$ such that $P(x(t),y(t),\dot{x}(t),\dot{y}(t))$
is Laurent series in $t$. This problem needs to be studied and
understood.
\end{Exmp}

\end{document}